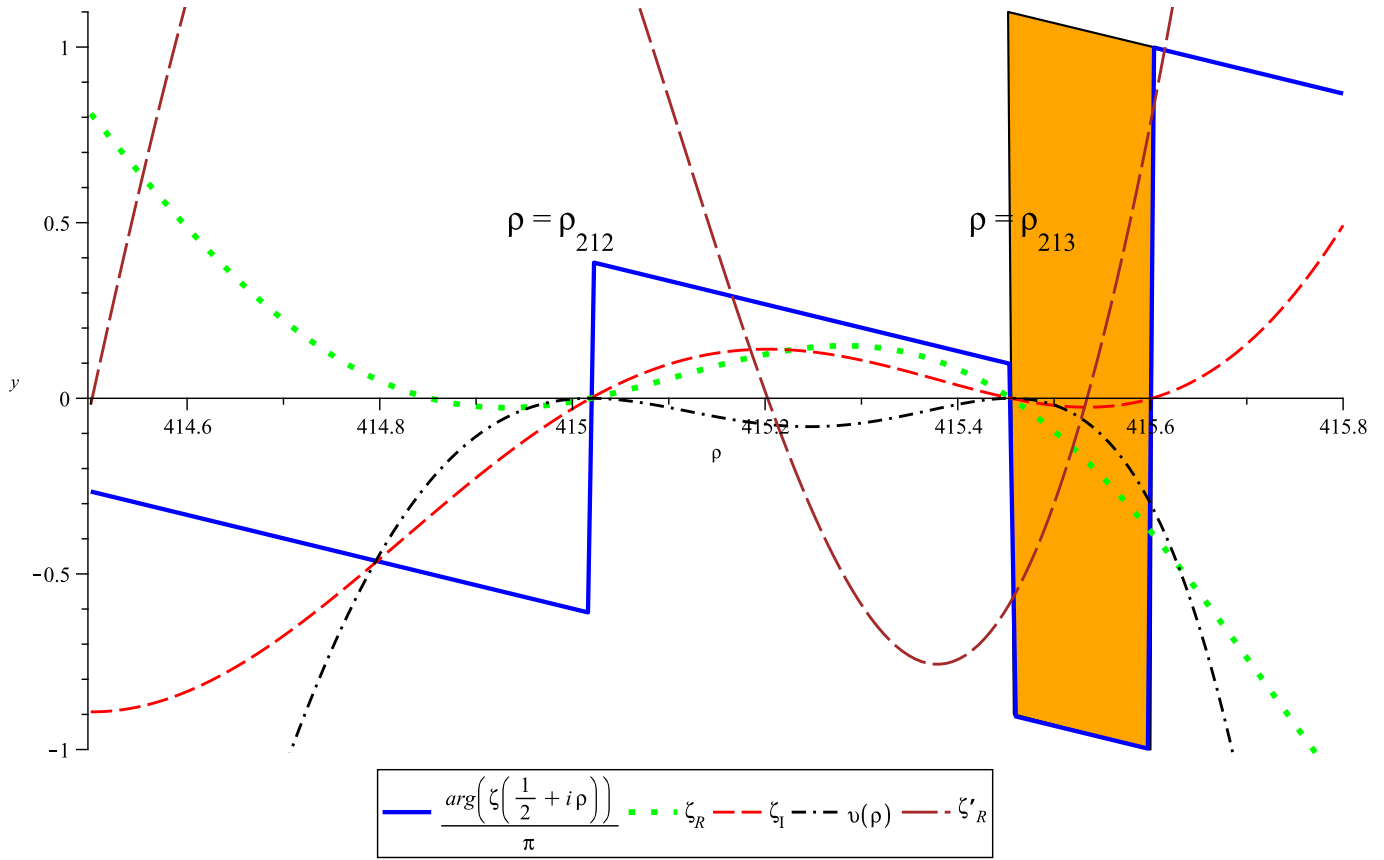

# Exploring Riemann's functional equation

Michael Milgram

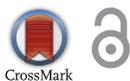

**cogent** · mathematics

**PURE MATHEMATICS | RESEARCH ARTICLE**

# Exploring Riemann's functional equation

Michael Milgram[1]*

**Abstract:** An equivalent, but variant form of Riemann's functional equation is explored, and several discoveries are made. Properties of Riemann's zeta function $\zeta(s)$, from which a necessary and sufficient condition for the existence of zeros in the critical strip, are deduced. This in turn, by an indirect route, eventually produces a simple, solvable, differential equation for $arg(\zeta(s))$ on the critical line $s = 1/2 + i\rho$, the consequences of which are explored, and the "LogZeta" function is introduced. A singular linear transform between the real and imaginary components of $\zeta$ and $\zeta'$ on the critical line is derived, and an implicit relationship for locating a zero ($\rho = \rho_0$) on the critical line is found between the arguments of $\zeta(1/2 + i\rho)$ and $\zeta'(1/2 + i\rho)$. Notably, the Volchkov criterion, a Riemann Hypothesis (RH) equivalent, is analytically evaluated and verified to be half equivalent to RH, but RH is not proven. Numerical results are presented, some of which lead to the identification of *anomalous zeros*, whose existence in turn suggests that well-established, traditional derivations such as the Volchkov criterion and counting theorems require re-examination. It is proven that the derivative $\zeta'(1/2 + i\rho)$ will never vanish on the perforated critical line ($\rho \neq \rho_0$). Traditional asymptotic and counting results are obtained in an untraditional manner, yielding insight into the nature of $\zeta(1/2 + i\rho)$ as well as very accurate asymptotic estimates for distribution bounds and the density of zeros on the critical line.

**Subjects:** Advanced Mathematics; Analysis - Mathematics; Complex Variables; Integral Transforms & Equations; Mathematical Analysis; Mathematics & Statistics; Number Theory; Science; Special Functions

## ABOUT THE AUTHOR

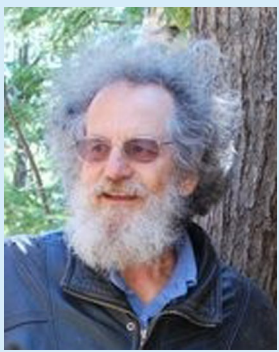

Michael Milgram

The author has spent most of his career studying and applying mathematics to problems in Physics and Engineering. This led him to develop some expertise in the mathematical field of Special Functions, which in turn led to an interest in the Zeta function, as a diversion from his everyday pursuit of solutions to problems in Monte Carlo simulation and computer-aided visualization of CNC machining.

## PUBLIC INTEREST STATEMENT

Prime numbers have fascinated mathematicians since the invention of arithmetic. One of the main discoveries by Bernhard Riemann, one of the greatest nineteenth-century mathematicians, was a link between what is called the "Riemann Zeta function" and the distribution of prime numbers. Since then, mathematicians (and others) have been studying zeta function properties and slowly unveiling its secrets. This is important because prime numbers form the basis for the major encryption algorithms in use today. If you indulge in Internet banking or commerce, you are using prime numbers.

In this paper, a new relationship is deduced that equates an arcane property of the zeta function to simpler mathematical entities whose properties are well known. As an anonymous referee has stated: "... such relation ... would lead to extremely important consequences". Additionally, a surprising numerical discovery is presented which calls into question many previously thought-to-be-well-founded results.

cogent ·· mathematics

## 1. Introduction

The course of other work has led to an exploration of a not-so-well-known variant of Riemann's functional equation. Thought to be equivalent to the classical equation that relates the functions $\zeta(s)$ and $\zeta(1-s)$, the variant studied here relates $\zeta(1-s)$ and its derivatives $\zeta'(s)$ and $\zeta'(1-s)$. As has been discovered elsewhere (e.g. Guillera, 2013, Eq. 10; Spira, 1973), this variant has proven to be a surprisingly rich source from which possible new properties of $\zeta(s)$ can be unearthed. At a minimum, the rediscovery of previously known results is obtained in a completely different manner from the usual textbook and literature approach. This report is a summary of those explorations.

In Section 2, the notation and other results drawn from the literature are summarized. Section 3 recalls the variant functional equation, from which it is possible to infer the existence of a functional relationship between $\zeta'(s)$ and $\zeta'(1-s)$ that yields necessary and sufficient implicit conditions for locating points $s = s_0$ such that $\zeta(s_0) = 0$. The main result of the next two sub-sections (4.1) and (4.2) is the derivation (and rediscovery) of two implicit, and later shown-to-be-equivalent, equations whose solution(s) locate zeros on the critical line $s = 1/2 + i\rho$. In a third sub-section (4.3), a singular linear transformation relating the real and imaginary components of $\zeta(1/2 + i\rho)$ and $\zeta'(1/2 + i\rho)$ is presented, generating a third implicit equation for locating the non-trivial zeros. Prior to presenting an analytic solution to the first two of these equations in Section 7, and noting the numerical equivalence of the third, the development digresses.

Section 5 introduces a related functional with interesting properties, and in Section 6, those properties are employed to obtain and solve a simple differential equation satisfied by $arg(\zeta(1/2 + i\rho))$ and thereby introduce the LogZeta function, in complete analogy to the well-known LogGamma function (England, Bradford, Davenport, & Wilson, 2013; Weisstein, 2005). Amalgamating all those results in Section 7 yields requisite conditions for the existence of $\rho_0$ such that $\zeta(1/2 + i\rho_0) = 0$, reproducing similar results obtained elsewhere (de Reyna & Van de Lune, 2014). In Section 8, by analytic integration, the properties of $arg(\zeta(1/2 + i\rho))$ are used to ostensibly verify that the Volchkov Criterion (Volchkov, 1995), advertised as "equivalent" to the Riemann Hypothesis (RH), possesses one-half of that property. Although this does not lead to a proof of RH, it does lead to some insight about the usefulness of so-called "Riemann equivalences" (Conrey & Farmer, n.d.). However, the subsequent Section 9 shows that much of the foregoing, as well as well-known classical results, are in need of revision because of the peculiar numerical properties of $\zeta(1/2 + i\rho)$, from which the phenomenon of *anomalous zeros* is identified. Additionally, the results from Section 6 lead to a proof (in Section 10) that $\zeta'(1/2 + i\rho) \neq 0$ on the perforated critical line without recourse to RH, almost closing the last gap in an investigation initiated many years ago (e.g. Spira, 1973). Section 11 investigates the asymptotic behaviour of $|\zeta(1/2 + i\rho)|$, again reproducing known results. In Section 12, these results are reused to deduce estimates that provide upper and lower bounds for locating the k'th zero on the critical line, as well as the density and maximum separation of such points. In a numerical diversion, the claimed location of zero number "googol" (França & LeClair (2013)) is tested; it is found that the value as declared should probably be googol-1, although its position on the critical line appears to be accurate within the resolution provided. (Aside: The word "googol" traditionally refers to the number $10^{100}$ and is easily confused with the similarly sounding word "google" usurped by a well-known search engine.)

Throughout, some derivations require considerable analytic perseverance to deal with expressions involving many terms. Such derivations are noted, and are left as exercises for the reader with the suggestion that a computer algebra program be utilized. In general, these computations, although lengthy, do not employ anything other than the use of well-known trigonometric and other identities involving Gamma and related functions (Olver, Lozier, Boisvert, & Clark, 2010). An example of such a calculation, in which the differential Equation (6.1) can be found as part of the supplemental material which may be accessed by a hyperlink following Section 14. A list of notations and symbols will be found in Appendix 1. The Maple computer code (Maplesoft, A Division of Waterloo Maple, 2014) is the source used for many of the calculations contained here. When that occurs, references to Maple are simply indicated in the text by the word Maple in parenthesis, in order that calculations can be reproduced,

## 2. Preamble

A number of known results and notation, required in forthcoming Sections, are quoted here. A basic result used throughout employs the polar form of $\zeta(s)$ as

$$\zeta(s) = e^{i\alpha(s)}|\zeta(s)| \tag{2.1}$$

differing by a sign from that used elsewhere (de Reyna & Van de Lune, 2014). Specific to the "critical line" $s = 1/2 + i\rho$, the specialized form is written

$$\zeta(1/2 + i\rho) = e^{i\alpha}|\zeta| = e^{i\alpha}\sqrt{\zeta_R^2 + \zeta_I^2} \tag{2.2}$$

where explicit dependence on the variable $\rho$ is usually omitted for clarity. Throughout, subscripts "R" and "I" refer to the real and imaginary parts of the associated symbol, respectively; all derivatives are taken with respect to $\rho$ and, always, $\rho \geq 0$. On the critical line, the arguments $\theta$ and $\beta$ of associated functions are defined by

$$\Gamma(1/2 + i\rho) = e^{i\theta}|\Gamma| \tag{2.3}$$

$$\zeta'(1/2 + i\rho) = e^{i\beta}|\zeta'| \tag{2.4}$$

For use in later sections, define

$$f(\rho) = \frac{4\cosh(\pi\rho)}{2\ln(2\pi)\cosh(\pi\rho) - 2\,\Re(\psi(1/2 + i\rho))\cosh(\pi\rho) + \pi}\,. \tag{2.5}$$

This function has a simple pole at $\rho_s \equiv \rho = 0.628...$, and for $\rho > \rho_s$, $f(\rho) < 0$. See Section 3.

Although the real functions $|\zeta|^2$ and $\alpha$ are (almost) independent, known relationships (de Reyna & Van de Lune, 2014, Proposition 7 and equivalently Milgram, 2011, Eqs. 3.1, 3.2 and Appendix B) exist between the real and imaginary components of $\zeta(1/2 + i\rho)$, specifically

$$\tan(\alpha) \equiv \frac{\zeta_I}{\zeta_R} = \mathfrak{P} \tag{2.6}$$

where

$$\mathfrak{P} = \frac{C_p\cos(\rho_\pi) + C_m\sin(\rho_\pi) - \sqrt{\pi}}{C_m\cos(\rho_\pi) - C_p\sin(\rho_\pi)}\,. \tag{2.7}$$

An equivalent form of (Equation 2.7) can be obtained by writing its various terms in polar form, using the representation (de Reyna & Van de Lune, 2014, Eq. 4)

$$\Gamma(1/2 + i\rho) = \sqrt{\frac{\pi}{\cosh(\pi\rho)}}\exp(i\,\theta)\,, \tag{2.8}$$

alternatively giving

$$\tan(\alpha) \equiv \frac{\zeta_I}{\zeta_R} = \frac{-\cosh(\frac{\pi\rho}{2})\cos(\rho_\theta) - \sinh(\frac{\pi\rho}{2})\sin(\rho_\theta)}{\sinh(\frac{\pi\rho}{2})\cos(\rho_\theta) - \cosh(\frac{\pi\rho}{2})\sin(\rho_\theta)}$$
$$+ \frac{\sqrt{\cosh(\pi\rho)}}{\sinh(\frac{\pi\rho}{2})\cos(\rho_\theta) - \cosh(\frac{\pi\rho}{2})\sin(\rho_\theta)}\,, \tag{2.9}$$

where the various symbols in the above are defined in Appendix 1. Curiously, relevant to the well-known result (Titchmarsh & Heath-Brown, 1986, Eq. 4.17.2)

$$arg(\zeta(1/2 + i\rho)) = arg(\Gamma(1/4 + i\rho/2)) - \frac{\rho}{2}\log(\pi) - k\pi\,,$$  (2.10)

and left as an exercise for the reader, it may also be shown that

$$\frac{\Im(\Gamma(1/4 + i\rho/2))}{\Re(\Gamma(1/4 + i\rho/2))} = -\mathfrak{P}$$  (2.11)

if $\rho_\pi \equiv \rho\log(2\pi)$ is replaced by $\rho\log(2)$ in the definition (Equation 2.7). As has been obtained elsewhere (Milgram, 2011, Eq. 3.7), in the case that $\zeta = \zeta_R = \zeta_I = 0$, corresponding to a non-trivial zero $\zeta(1/2 + i\rho_0) = 0$, the ratio of those quantities appearing on the left-hand side of both (Equation 2.7 and 2.9) satisfies

$$\lim_{\rho \to \rho_0} \frac{\zeta_I}{\zeta_R} = p\left(\frac{\zeta^{(n)}_I}{\zeta^{(n)}_R}\right)^\rho$$  (2.12)

when the zero is of order $n$, where $p = (-1)^n$, in terms of the appropriate $n^{th}$ derivatives, by invoking l'Hôpital's rule. Also see Equations (6.3 and 6.4) below. Another result that will prove useful is obtained by direct differentiation of Equation (2.2), that being

$$\frac{|\zeta'(1/2 + i\rho)|}{|\zeta(1/2 + i\rho)|} = \frac{\alpha'}{\cos(A - B)}\,,$$  (2.13)

where

$$A \equiv arg(\zeta(1/2 + i\rho))$$
$$B \equiv arg(\zeta'(1/2 + i\rho))\,.$$  (2.14)

An important distinction must be made between discontinuous functions corresponding to, and denoted by, the notation "arg" and their continuous, multi-sheeted counterparts—the two entities differ by a constant equal to $k\pi$ and $k$ is always an integer. As an example, consider Figure 1 where the imaginary part of the multi-sheeted function (Weisstein, 2005) LogGamma$(1/2 + i\rho)$ is compared to the discontinuous function $arg(\Gamma(1/2 + i\rho))$, in both cases as a function of $\rho$. The difference, normalized by $\pi$ is also shown, to demonstrate that it is always an integer "constant", which, however, changes at each discontinuous point of the $arg$ operator. As indicated in the Figure, discontinuities of the argument are associated with a sign change of $\Im(\Gamma(1/2 + i\rho))$ coincident with $\Re(\Gamma(1/2 + i\rho)) < 0$—that is, whenever $arg(\Gamma(1/2 + i\rho)) = \pm\pi$, and the discontinuity is $2\pi$. In Section 9, a similar comparison is made using the function $\zeta(1/2 + i\rho)$.

**Figure 1.** Comparison between $\Im(\text{LogGamma}(1/2 + i\rho))$ and $arg(\Gamma(1/2 + i\rho))$, along with the normalized associated real and imaginary parts of $\Gamma(1/2 + i\rho)$.

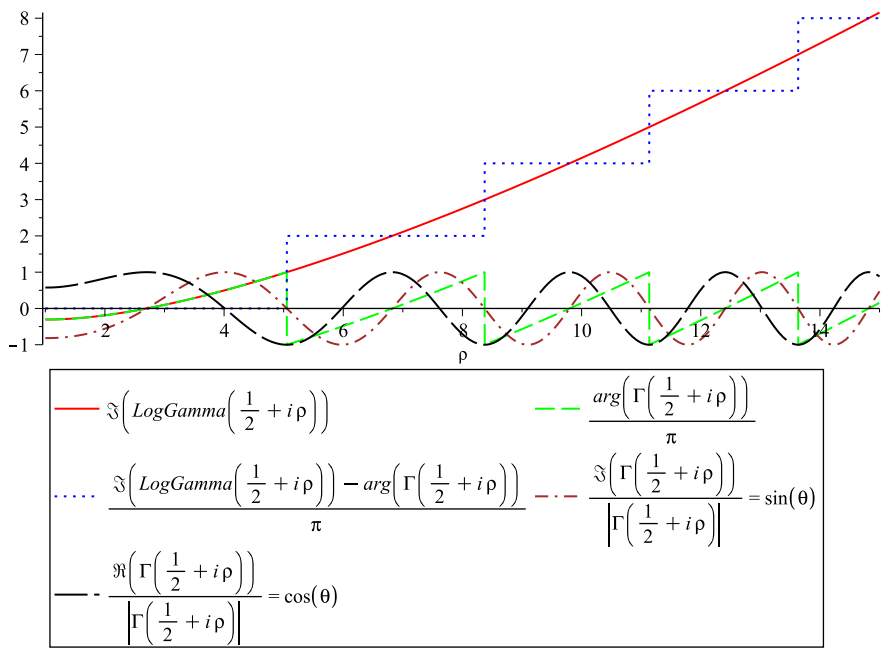

### 3. The basic functional equation

The following proofs of sufficiency and necessity are valid if $\zeta(s_0) = 0$ is a simple zero. See Section 10.

#### 3.1. A sufficient condition that $\zeta(s) = 0$

From the functional equation for $\zeta(s)$, that is,

$$\zeta(1 - s) = \frac{2\,\Gamma(s)\,\cos(\pi s/2)\,\zeta(s)}{(2\pi)^s} \tag{3.1}$$

we have (Guillera, 2013, Eq. 10 - misprinted; Spira, 1973, Eq. 1)

$$\zeta'(1 - s) + \zeta'(s)\frac{2\,\Gamma(s)\cos(\pi s/2)}{(2\,\pi)^s} = (\ln(2\,\pi) - \psi(s) + \frac{\pi}{2}\tan(\pi s/2))\zeta(1 - s) \tag{3.2}$$

valid for all $s$. Alternatively, in equivalent form, define the normalized right- and left-hand sides of Equation (3.2) by

$$\mathfrak{T}(s) \equiv (\ln(2\,\pi) - \psi(s) + \frac{\pi}{2}\tan(\pi s/2))\,\zeta(1 - s)/\zeta'(s) \tag{3.3}$$

and

$$\mathfrak{L}(s) \equiv \frac{\zeta'(1 - s)}{\zeta'(s)} + 2\cos(\pi s/2)\,\Gamma(s)\,(2\,\pi)^{(-s)}, \tag{3.4}$$

in which case Equation (3.2) can be written

$$\mathfrak{L}(s) = \mathfrak{T}(s). \tag{3.5}$$

Since it is known (Spira, 1973) that $\zeta'(s) \neq 0$ in the open critical half-strip $0 \leq \Re(s) < 1/2$, the following requires that $s$ be constrained to that region, although the known symmetry imposed by the functional equation, viz. $\zeta(1 - s_0) = 0$ implies $\zeta(s_0) = 0$ and the reverse, means that the following can be generalized to $\Re(s) > 1/2$. We now show that $\mathfrak{L}(s_0) = 0$ is a necessary and sufficient condition for $\zeta(s_0) = 0$. Suppose $s = s_0$ such that

$$\frac{\zeta'(1 - s_0)}{\zeta'(s_0)} = -2\cos(\pi s_0/2)\,\Gamma(s_0)\,(2\pi)^{-s_0}. \tag{3.6}$$

Clearly, if Equation (3.6) is true, or, alternatively

$$\mathfrak{L}(s_0) = 0, \tag{3.7}$$

then the right-hand side of Equation (3.2) vanishes; that is,

$$\zeta(1 - s_0) = 0 \tag{3.8}$$

unless, for that same value of $s_0$, the factor

$$\ln(2\pi) - \psi(s_0) + \frac{1}{2}\pi\tan(\frac{\pi s_0}{2}) = 0 \ . \tag{3.9}$$

Is it possible that Equations (3.7 and 3.9) can be simultaneously true for some value(s) of $s_0$? On the perforated critical line, where it will shortly be shown (see Section 10) that $\zeta'(1/2 + i\rho) \neq 0$, solving Equation (3.9) requires that $\rho$ must satisfy

$$\Re(\psi(1/2 + i\,\rho)) = \frac{\ln(2\pi)\cosh(\pi\,\rho) + \pi/2}{\cosh(\pi\,\rho)} \tag{3.10}$$

and from Figure 2, we see that this condition is only satisfied at a single point $\rho_s = 6.2898...$ because the right-hand side of Equation (3.10) approaches $\ln(2\pi)$ asymptotically, whereas the left-hand side is monotonically increasing as $\log(\rho)$ (Olver et al., 2010, Eq. 5.11.2). In this case, the imaginary counterpart of Equation (3.9) vanishes. Because this point does not coincide with a zero of $\zeta(1/2 + i\rho)$, on the critical line, Equation (3.6), or equivalently Equation (3.7), is a sufficient condition for Equation (3.8) to be true, unless $\rho_0 = 6.2898...$, a constant quoted to many significant figures in de Reyna and Van de Lune (2014, Corollary 9). See also Section 10 - Note 1.

Now, consider the general case $s = \sigma + i\rho$ with $0 < \sigma \leq 1/2$—the lower half of the "critical strip". A hypothetical solution of Equation (3.9) requires the existence of solutions $(\sigma_0, \rho_0)$ simultaneously satisfying

$$\Re(\psi(\sigma_0 + i\,\rho_0)) = \ln(2\pi) + \frac{1}{2}\frac{\pi\sin(\pi\,\sigma_0)}{\cos(\pi\,\sigma_0) + \cosh(\pi\,\rho_0)} \tag{3.11}$$

and

$$\Im(\psi(\sigma_0 + i\,\rho_0)) = \frac{1}{2}\frac{\pi\sinh(\pi\,\rho_0)}{\cos(\pi\,\sigma_0) + \cosh(\pi\,\rho_0)} \ . \tag{3.12}$$

**Figure 2. Numerical demonstration of Equation (3.10).**

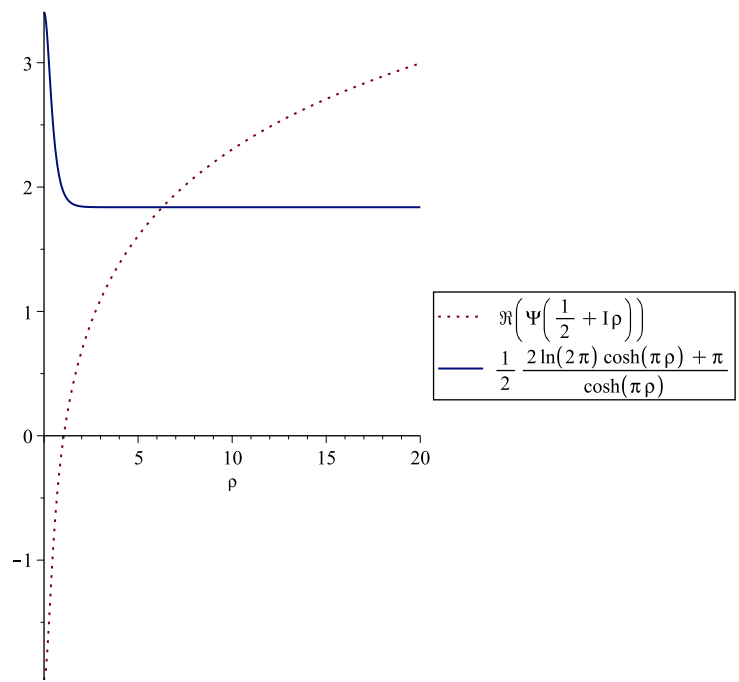

**Figure 3.** Numerical demonstration of Equation (3.4) for two values of $\sigma$ near the first two zeros of $\zeta(1/2 + i\rho)$.

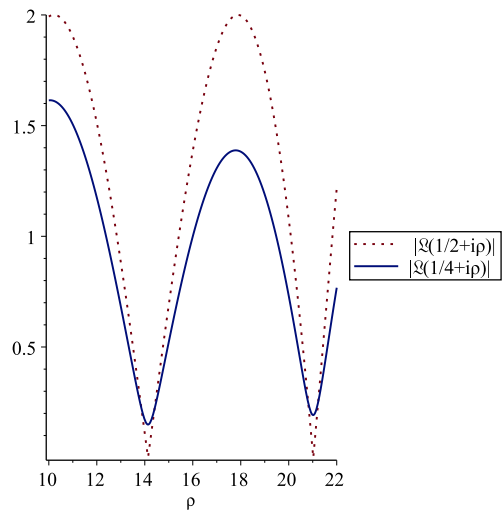

If $\sigma_0 = 1/2$, Equation (3.11) reduces to Equation (3.10) and (3.12) reduces to an identity (NIST Digital Library of Mathematical Functions, 2014, Eq. 5.4.17). A numerical study of Equation (3.11) for $0 < \sigma < 1/2$ shows that a single solution exists when $\rho_0 \approx 6.2$, similar to that shown in Figure 2, but otherwise, there are no numerical solutions for the previous reasons. Similarly, the left-hand side of Equation (3.12) asymptotically approaches its right-hand side as $\rho_0 \to \infty$, but the equality is never satisfied. Therefore, the existence of a point $s_0$ solving Equation (3.7) is a sufficient condition for a zero of $\zeta(s)$ to exist because the Equations (3.7 and 3.9) do not simultaneously vanish.

By a way of illustration, consider Figure 3 which shows $|\mathfrak{L}(\sigma + i\rho)|^2$ for two different choices of $\sigma$ near the first two known non-trivial zeros of $\zeta(1/2 + i\rho)$. Because of the above, this Figure immediately suggests that $\mathfrak{L}(\sigma + i\rho) \neq 0$ unless $\sigma = 1/2$.

### 3.2. A necessary condition that $\zeta(s) = 0$

From the functional Equation (3.1), $\zeta(1 - s_0) = 0$ implies $\zeta(s_0) = 0$. If $\zeta(1 - s_0) = 0$, then $\mathfrak{T}(s_0) = 0$, since the other factors of $\mathfrak{T}(s)$ do not diverge if $\mathfrak{F}(s) \neq 0$ (corresponding to the existence of the so-called trivial zeros). From Equation (3.5), this implies $\mathfrak{L}(s_0) = 0$ and necessity follows because $\zeta'(s) \neq 0$ except possibly at a zero (see Section 10). This proves necessity subject to the simplicity of a zero.

## 4. Implicit specification of zeros along the critical line

### 4.1. First implicit specification of the zeros along the critical line

Inspired by Figure (3) and Equation (3.7) that together demonstrate that implicit numerical solutions of the equation $|\mathfrak{L}(1/2 + i\rho)|^2 = 0$ locate the zeros of $\zeta(1/2 + i\rho)$, consider the related equation

$$|\mathfrak{L}(1/2 + i\rho)|^2 = |\mathfrak{T}(1/2 + i\rho)|^2 \tag{4.1}$$

After considerable effort (Maple) to simplify (Equation 4.1) (left as an exercise for the reader with access to a computerized algebraic manipulation program), we find

$$\frac{4\,|\zeta|^2}{f(\rho)^2\,|\zeta'|^2} = 2 + \frac{2\left(\sinh(\frac{\pi\rho}{2})\sin(2\,\beta + \rho_\theta) + \cosh(\frac{\pi\rho}{2})\cos(2\,\beta + \rho_\theta)\right)}{\sqrt{\cosh(\pi\rho)}}\,. \tag{4.2}$$

For any value of $\rho = \rho_0$ corresponding to $|\zeta|^2 = 0$ and $|\zeta'|^2 \neq 0$ (see Section 10) on the left-hand side of Equation (4.2), the right-hand side of this equation must accordingly vanish as well because the remaining factors on the left-hand side never vanish if $\rho > \rho_s$. This produces an implicit equation that locates $\zeta(1/2 + i\rho_0) = 0$:

$$\frac{\sinh(\pi \rho/2)\sin(2\beta + \rho_\theta) + \cosh(\pi \rho/2)\cos(2\beta + \rho_\theta)}{\sqrt{\cosh(\pi \rho)}} = -1 \, . \tag{4.3}$$

This will be discussed further in Section 7.

### 4.2. Second implicit specification of the zeros along the critical line

Again inspired by Figure (3), we obtain a second implicit equation whose zeros coincide with the non-trivial zeros of $\zeta(1/2 + i\rho)$. Working from the real and imaginary parts of $\mathfrak{L}(1/2 + i\rho)$

$$\mathfrak{L}_R = \frac{\sin(\rho_\theta)\sinh(\frac{\pi\rho}{2}) + \cos(\rho_\theta)\cosh(\frac{\pi\rho}{2})}{\sqrt{\cosh(\pi\rho)}} + \frac{{\zeta'_R}^2 - {\zeta'_I}^2}{|\zeta'|^2} \tag{4.4}$$

$$\mathfrak{L}_I = \frac{\sin(\rho_\theta)\cosh(\frac{\pi\rho}{2}) - \cos(\rho_\theta)\sinh(\frac{\pi\rho}{2})}{\sqrt{\cosh(\pi\rho)}} - \frac{2\,\zeta'_R\,\zeta'_I}{|\zeta'|^2} \tag{4.5}$$

consider the function

$$\mathfrak{L}_1(\rho) \equiv \frac{d}{d\rho}|\mathfrak{L}(1/2 + i\rho)|^2 \tag{4.6}$$

at least some of whose zeros $\rho = \rho_0$ must coincide with the maxima and minima of $|\mathfrak{L}(1/2 + i\rho)|$.

Following a lengthy calculation (Maple), we find

$$\mathfrak{L}_1(\rho) = 2\,\sqrt{\frac{1}{\cosh(\pi\,\rho)}}\left[-\frac{1}{f(\rho)}|\zeta'|^2 + (\sin(\beta)\,\zeta''_I + \cos(\beta)\,\zeta''_R)\,|\zeta'|\right]T_1(\rho) \tag{4.7}$$

where again the various symbols in Equation (4.7) are defined in Appendix 1.

Figure 4 demonstrates that $\mathfrak{L}_1(\rho)$ passes through the first two zeros of $\zeta(1/2 + i\rho)$ with positive slope, and, as expected, alternating solutions corresponding to zeros of Equation (4.7) implicitly define those points $\rho = \rho_0$ such that $\zeta(1/2 + i\rho_0)$=0, consistent with the assumption that $|\mathfrak{L}(1/2 + i\rho)|^2$ does not contain intermediate maxima, minima or inflections. Of particular interest is the fact that Equation (4.7) consists of two factors that could potentially vanish. From a numerical study, the first of these, enclosed in square brackets ([...]), appears to be negative for all values of $\rho$.

**Figure 4. Numerical demonstration of Equation (4.8) near the first two zeros of $\zeta(1/2 + i\rho)$.**

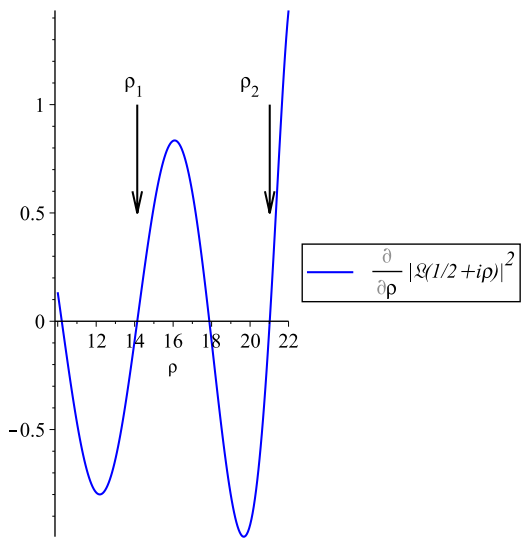

cogent •• mathematics

**Figure 5.** Numerical demonstration of Equation (4.8) near the first few zeros of $\zeta(1/2 + i\rho)$.

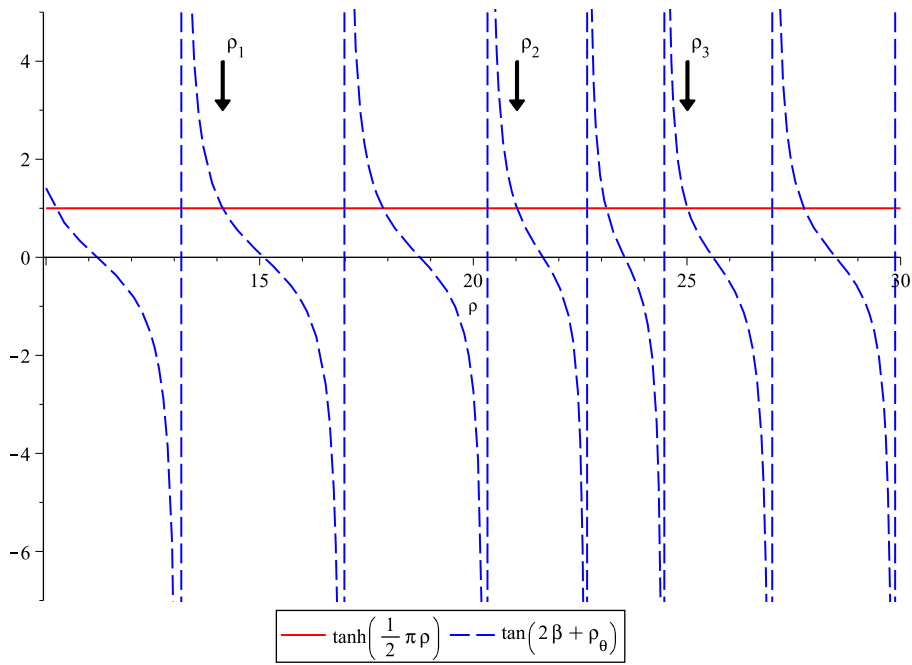

If so, the factor $T_1(\rho)$ (see Appendix 1) carries all the zeros of $\mathfrak{L}_1(\rho)$, from which we conclude that solutions of $\zeta(1/2 + i\rho_0) = 0$ correspond to alternating solutions of

$$\tanh\left(\frac{\pi\rho}{2}\right) = \tan\left(2\beta + \rho_\theta\right). \tag{4.8}$$

See Figure 5. This will be discussed further in Section 7.

### 4.3. Third implicit specification of the zeros along the critical line
With $s = 1/2 + i\rho$, from Equation (3.2 and 2.5), we find

$$\zeta_R/f(\rho) = a(\rho)\,\zeta'_R + b(\rho)\,\zeta'_I \tag{4.9a}$$

$$\zeta_I/f(\rho) = b(\rho)\,\zeta'_R + c(\rho)\,\zeta'_I \tag{4.9b}$$

where

$$a(\rho) = -\frac{1}{2}\cos(\rho\ln(\pi)) + \frac{1}{2} + \frac{\mathcal{F}_I\mathcal{F}_R\sin(\rho\ln(\pi)) + \mathcal{F}_R^2\cos(\rho\ln(\pi))}{\left|\Gamma(\frac{1}{4} + \frac{1}{2}i\rho)\right|^2} \tag{4.10}$$

$$b(\rho) = -\frac{1}{2}\sin(\rho\ln(\pi)) + \frac{-\mathcal{F}_I\mathcal{F}_R\cos(\rho\ln(\pi)) + \mathcal{F}_R^2\sin(\rho\ln(\pi))}{\left|\Gamma(\frac{1}{4} + \frac{1}{2}i\rho)\right|^2} \tag{4.11}$$

$$c(\rho) = 1 - a(\rho) \tag{4.12}$$

and we reiterate (from Appendix 1) that $\mathcal{F} \equiv \Gamma(1/4 + i\rho/2)$.

It is easily shown that the determinant of the transformation matrix Equation (4.9) is identically zero, and therefore the transformation is singular. This suggests that in some sense, the components of the derivative function $\zeta'(1/2 + i\rho)$ are more fundamental than are the components of the

cogent • mathematics

primary function $\zeta(1/2 + i\rho)$. That is, Riemann's functional equation in the form Equation (3.2) defines the latter from the former, but not the other way around, suggesting that something has been lost when Equation (3.2) is generated from the usual form of the functional Equation (3.1). Why should this be? This may be analogous to the fact that the act of differentiation followed by integration always introduces an arbitrariness in the form of a constant. Or, it may be due to the fact that differentiating (Milgram, 2011, Eqs. 3.1 and 3.2) with respect to $\rho$ will always give a result where $\zeta'_R$ and $\zeta'_I$ are a mixture of $\zeta_R$ and $\zeta_I$ that cannot be disentangled in the form of an inverse transform of Equations (4.9a) and (4.9b).

In Section 10, it is shown that $\zeta'(1/2 + i\rho) \neq 0$ with conditions. Thus, setting the left-hand sides of Equations (4.9a and 4.9b) to zero will result in a relationship between the various components of these equations at a zero of $\zeta(1/2 + i\rho)$. Specifically, $\rho = \rho_0$ satisfying the following

$$\sin(\rho_0 \ln(\pi)) = \frac{2\,(\mathcal{F}_I\,\zeta'_R + \mathcal{F}_R\,\zeta'_I)\,(\mathcal{F}_I\,\zeta'_I - \mathcal{F}_R\,\zeta'_R)}{|\zeta'|^2\,|\mathcal{F}|^2} \tag{4.13}$$

and

$$\cos(\rho_0 \ln(\pi)) = \frac{4\,\mathcal{F}_I\,\mathcal{F}_R\,\zeta'_I\,\zeta'_R + (\zeta'^2_I - \zeta'^2_R)(-\mathcal{F}^2_I + \mathcal{F}^2_R)}{|\zeta'|^2\,|\mathcal{F}|^2} \tag{4.14}$$

defines all the non-trivial zeros on the critical line plus one extra zero associated with the factor $1/f(\rho_0)$. The results (Equations 4.13 and 4.14) can be consolidated by equating the ratio of the two left-hand sides in the form

$$\arctan(\sin(\rho_0 \log(\pi)), \cos(\rho_0 \log(\pi))) \tag{4.15}$$

with a similar expression for the right-hand sides. The standard two parameter $\arctan(y, x)$ function is specified to consistently account for the signs of the various terms. Alternatively, solving Equations (4.13 and 4.14) for the ratio $\zeta'_I/\zeta'_R$ gives

$$\tan(\beta_0) \equiv \frac{\zeta'_I}{\zeta'_R} = \frac{|\mathcal{F}|^2\cos(\rho_0\ln(\pi)) - \mathcal{F}^2_I + \mathcal{F}^2_R}{-\sin(\rho_0\ln(\pi))\,|\mathcal{F}|^2 + 2\,\mathcal{F}_I\,\mathcal{F}_R} \tag{4.16}$$

corresponding to a solution for $\beta_0 = arg(\zeta'(1/2 + i\rho_0))$ at a non-trivial zero $\rho = \rho_0$. See Section 7.

## 5. A related functional
Formally replace $\zeta'$ with $\zeta$ in Equation (3.4) and using Equation (3.1), define and evaluate

$$\mathfrak{M}(s) \equiv \frac{\zeta(1-s)}{\zeta(s)} + 2\cos(\pi s/2)\Gamma(s)(2\pi)^{(-s)} = 4\cos(\pi s/2)\Gamma(s)(2\pi)^{(-s)}. \tag{5.1}$$

For the particular case $s = 1/2 + i\rho$, using standard trigonometric and Gamma function identities (e.g. Olver et al., 2010, Eq. 5.4.4), it is easy to establish that

$$|\mathfrak{M}(1/2 + i\rho)|^2 = 4. \tag{5.2}$$

It now becomes possible to use (Equation 5.2) to discover interesting results because of Equation (3.5). Define the modified right-hand side of Equation (3.5) corresponding to Equation (5.1) with $s = 1/2 + i\rho$ by adding and subtracting the terms that convert $\mathfrak{L}(s)$ into $\mathfrak{M}(s)$, specifically

$$\mathfrak{Q}(1/2 + i\rho) \equiv \mathfrak{T}(1/2 + i\rho) + \frac{\zeta(1/2 - i\,\rho)}{\zeta(1/2 + i\,\rho)} - \frac{\zeta'(1/2 - i\,\rho)}{\zeta'(1/2 + i\,\rho)}. \tag{5.3}$$

Then, from Equation (5.2), we have

$$|\mathfrak{Q}(1/2 + i\rho)|^2 = 4 \ .$$

(5.4)

## 6. $arg(\zeta(1/2 + i\rho))$ **and related entities**

After algebraic simplification (Maple), Equation (5.4) eventually yields

$$|\zeta|^2 = \zeta_I^2 + \zeta_R^2 = (\zeta'_I \, \zeta_I + \zeta'_R \, \zeta_R) f(\rho) \ .$$

(6.1)

(Background note: This is a fairly lengthy, but straightforward calculation. First, break $\mathfrak{Q}$ into its real and imaginary parts using Equations (3.3 and 5.3), evaluate $|\mathfrak{Q}|^2 = \mathfrak{Q}_R^2 + \mathfrak{Q}_I^2 = 4$ and simplify the resulting Equation (5.4). The simplification sequence involves nothing more than the application of well-known identities involving $\Gamma$ and related functions. Due to the complexity and length of this calculation, the details are left as an exercise for the reader, who may consult a supplementary file which can be accessed by a hyperlink following Section 14—that being an annotated Maple worksheet used to reproduce the calculation. A second, equally lengthy, but independent derivation of Equation (6.1) has recently been obtained—ms. in preparation).

Setting $\rho = 0$ in Equation (6.1) reproduces a known relationship between $\zeta(1/2)$ and $\zeta'(1/2)$. Furthermore, noticing that

$$\left(\frac{\zeta_I}{\zeta_R}\right)' = \frac{(\zeta_I)'}{\zeta_R} - \frac{\zeta_I \, (\zeta_R)'}{\zeta_R^2}$$

(6.2)

by recalling that

$$(\zeta_I)' = (\zeta')_R \equiv \zeta'_R$$

(6.3)

and

$$(\zeta_R)' = -(\zeta')_I \equiv -\zeta'_I$$

(6.4)

we recognize Equation (6.1) to be a simple differential equation

$$\frac{d\,g(\rho)}{d\rho} = \frac{1 + g(\rho)^2}{f(\rho)}$$

(6.5)

or, equivalently

$$\alpha'(\rho) \equiv \frac{d\,\alpha(\rho)}{d\,\rho} = 1/f(\rho)$$

(6.6)

where

$$g(\rho) \equiv \frac{\zeta_I}{\zeta_R} = \tan(\alpha(\rho)).$$

(6.7)

Integrating (Equation 6.6) between $\rho = \rho_1$ and $\rho = \rho_2$ gives

$$arg(\zeta(1/2 + i\,\rho_2)) - arg(\zeta(1/2 + i\,\rho_1)) = \frac{1}{2}(\rho_2 - \rho_1)\ln(2\,\pi)$$

$$-\frac{1}{2}\int_{\rho_1}^{\rho_2} \mathfrak{R}(\psi(1/2 + i\,\rho))\,d\rho - \frac{1}{2}\,\arctan(e^{\pi\,\rho_1}) + \frac{1}{2}\,\arctan(e^{\pi\,\rho_2}) - (k+2)\,\pi$$

(6.8)

where the arbitrary constant $(k + 2)\pi$ has been chosen such that $k = 0$ corresponds to a consistent answer when $\rho_1, \rho_2 \approx 0$. This result is equivalent to de Reyna and Van de Lune (2014, Eq. 8) when $\rho_1 = 0$ and $\rho_2 = \rho$, in which case we find

$$\arg(\zeta(1/2 + i\rho)) = -\frac{1}{2}\int_0^\rho \Re(\psi(1/2 + it))\,dt + \frac{\rho}{2}\ln(2\pi) - \frac{9\pi}{8} + \frac{1}{2}\arctan(e^{\pi\rho}) + k\pi. \qquad (6.9)$$

It is important to recognize that $\alpha(\rho)$ in Equation (6.6) is a continuous function, whereas its companion $arg(\zeta(1/2 + i\rho))$ is discontinuous, and therefore the solution (Equations 6.8 or 6.9) only applies in a multi-sheeted sense. That is, in order to satisfy Equation (6.6), in each of Equation (6.8 or 6.9), $k$ serves as a *local* variable, constant over a range where $arg(\zeta(1/2 + i\rho))$ is continuous, but effectively, $k = k(\rho)$ globally. See Section 9 for further discussion, and Section 8 for an application of this point.

The integral in Equation (6.9) can be evaluated analytically, using the well-known expansion (NIST Digital Library of Mathematical Functions, 2014, Eq. 5.7.6)

$$\psi(1/2 + x) = -\gamma - \frac{1}{1/2 + x} - \sum_{n=1}^\infty \left(\frac{1}{n + 1/2 + x} - \frac{1}{n}\right). \qquad (6.10)$$

First, temporarily omit the $\Re(...)$ operator, then interchange the sum and integration operators, and, after the (trivial) integration is accomplished, compute the real part of the result, and sum the resulting series utilizing (Hansen, 1975, Eq. 42.1.5)

$$\sum_{n=0}^\infty \left(\arctan(\frac{y}{n + x}) - \frac{y}{n + x}\right) = y\,\psi(x) - \arg(\Gamma(x + iy)),\, x + iy \neq 0, -1, -2 \dots. \qquad (6.11)$$

to obtain (see de Reyna & Van de Lune, 2014, Section 3; also see Section 12, Equation 12.2, below where it is pointed out that this result can be simply obtained by integrating by parts)

$$\int_0^\rho \Re\left(\psi\left(\frac{1}{2} + it\right)\right)dt = \Im(\,\mathrm{Log}\Gamma(1/2 + i\rho)) = \arg\left(\Gamma\left(\frac{1}{2} + i\rho\right)\right) + k\pi. \qquad (6.12)$$

(Digression: It is worth noting that an equivalent, but unevaluated, form of the sum in Equation (6.11) arises in de Reyna and Van de Lune (2014, Eq. 9) where it is suggested that it is perhaps new. In Hansen (1975), this result is attributed to Abramowitz and Stegun (1964, Eq. 6.1.27), and, somewhat surprisingly, it appears to have been omitted from (NIST Digital Library of Mathematical Functions, 2014). See also (Weisstein, 2005) where this sum is used to define the LogGamma function. Of course

$$arg(\Gamma(1/2 + i\rho)) = \Im(\log(\Gamma(1/2 + i\rho))) \qquad (6.13)$$

and the multi-sheeted LogGamma function of argument $(1/2 + i\rho)$ differs from the function $\log(\Gamma(1/2 + i\rho))$ by an additive term equal to $k\pi$.) With this result, Equation (6.8) now becomes

$$\arg(\zeta(\frac{1}{2} + i\rho_2)) - \arg(\zeta(\frac{1}{2} + i\rho_1)) = \frac{1}{2}(\rho_2 - \rho_1)\ln(2\pi) + \frac{1}{2}\arg(\Gamma(\frac{1}{2} + i\rho_1))$$
$$-\frac{1}{2}\arg(\Gamma(\frac{1}{2} + i\rho_2)) - \frac{1}{2}\arctan(e^{\pi\rho_1}) + \frac{1}{2}\arctan(e^{\pi\rho_2}) + k\pi \qquad (6.14)$$

and, in the case that $\rho_1 = 0, \rho_2 = \rho$,

$$\arg(\zeta(\frac{1}{2} + i\rho)) = \frac{\rho}{2}\ln(2\pi) - \frac{1}{2}\arg(\Gamma(\frac{1}{2} + i\rho)) - \frac{9\pi}{8} + \frac{1}{2}\arctan(e^{\pi\rho}) - k\pi. \qquad (6.15)$$

Notice that $\alpha'(\rho)$ in Equation (6.6) and $\alpha(\rho)$ in Equation (6.14) differ from substitute symbols ($\alpha \to A$ and $\beta \to B$) in Equations (2.13 and 2.14) because of the ambiguity associated with the term $k\pi$. Also, the discontinuities on both sides of Equations (6.14 and 6.15) do not coincide—for any value of $\rho$, a corresponding value of $k$ must be carefully chosen. From Equation (6.9), it becomes convenient to define the multi-sheeted (i.e. continuous) LogZeta function (see Figure 6):

cogent •• mathematics

**Figure 6. Comparison of the right- and left-hand sides of Equation (6.9) with $k = 1$.**

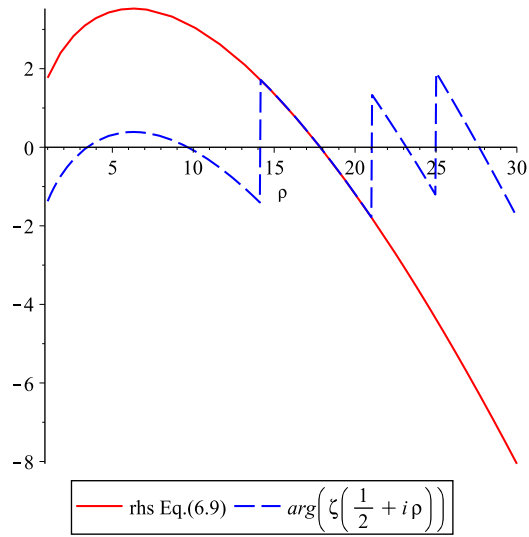

$$\text{LogZ}(1/2 + i\rho) \equiv -\log\Gamma(1/2 + i\rho)/2 + i(\tfrac{\rho}{2}\log(2\pi) - 9\pi/8 + \arctan(e^{\pi\rho})/2). \tag{6.16}$$

In this form, and throughout this work, $\alpha$ associates with LogZeta function via

$$\alpha(\rho) = \Im(\log Z(1/2 + i\rho)), \tag{6.17}$$

and $\theta$ represents the LogGamma function in the same way (England et al., 2013). Using this definition, we have

$$arg(\zeta(1/2 + i\rho)) = \Im(\log(\zeta(1/2 + i\rho))) = \Im(\text{LogZ}(1/2 + i\rho)) + k\pi. \tag{6.18}$$

Finally, from Equation (2.7), we obtain

$$\tan(\Im(\text{LogZ}(1/2 + i\rho))) = \mathfrak{P} \tag{6.19}$$

and an equivalent form if Equation (2.9) is used.

## 7. Locating the zeros on the critical line

Write Equation (4.3) in generic form

$$\frac{\sinh(\tfrac{\pi\rho}{2})\sin(Y) + \cosh(\tfrac{\pi\rho}{2})\cos(Y)}{\sqrt{\cosh(\pi\rho)}} = -1. \tag{7.1}$$

Solutions to Equation (7.1) exist when

$$Y = -\frac{(2n+1)\pi}{4} + \kappa\pi + \arctan(e^{\pi\rho}) \tag{7.2}$$

provided that (when n is even)

$$(-1)^{n/2}\cos(\kappa\pi) = -1 \tag{7.3}$$

or (when n is odd)

$$(-1)^{(n-1)/2} \sin(\kappa\,\pi) = -1\,.  \tag{7.4}$$

By comparison of Equation (7.2) with Equations (4.3 and 6.15), we identify

$$\kappa = 2(\beta - \alpha)/\pi + n/2 + 2(k-1)\,,  \tag{7.5}$$

consistent with de Reyna and Van de Lune (2014, Proposition 17). Taking advantage of the arbitrariness of the *arg* operator modulo $2\pi$, it is thus possible to solve Equation (7.1) in a simpler form and conclude that alternating solutions of

$$\arg(\zeta'(\tfrac{1}{2} + i\,\rho)) - \arg(\zeta(\tfrac{1}{2} + i\,\rho)) = (K + 1/2)\pi\,,  \tag{7.6}$$

locate zeros of $\zeta(1/2 + i\rho)$ on the critical line, where $K$ is an integer, a result that is easily verified numerically, or analytically for any combination of $n$ and (constrained) $\kappa$ in Equations (7.3 and 7.4). Equation (7.6) also reproduces de Reyna and Van de Lune (2014, Proposition 13).

Insofar as alternative means of locating zeros are concerned, it is easily shown by simple trigonometric identities that Equation (7.2) is also a valid solution to Equation (4.8); therefore Equation (4.8) is algebraically equivalent to Equation (4.3). Although not obvious, it turns out that Equation (4.16) is also (numerically) equivalent to Equations (4.3 and 4.8). In contrast, solutions satisfying

$$\arg(\zeta'(\tfrac{1}{2} + i\,\rho)) - \arg(\zeta(\tfrac{1}{2} + i\,\rho)) = K\pi  \tag{7.7}$$

correspond to the in-between solutions of Equation (7.1), i.e. the maxima of $|\mathfrak{L}(1/2 + i\rho)|^2$ in Figure 3. Any of Equation (4.3), Equation (4.8) or Equation (4.16), treated as implicit equations, are numerically sensitive means of locating a zero because the two sides intersect one another at differing slopes (e.g. Figure 5), as opposed to treating either side of Equation (4.1) implicitly because those "intersections" occur at tangent points. Another effective numerical variation, also equivalent to Equation (2.12), yields the following simple corollary to Equation (7.6):

$\zeta(1/2 + i\rho_0) = 0$ whenever $\rho = \rho_0$ satisfies

$$1 + \tan(\alpha(\rho))\,\tan(\beta(\rho)) = 0.  \tag{7.8}$$

Any of these prescriptions, employed as a numerical means of locating zeros, give the appearance of requiring knowledge of both $\alpha$ and $\beta$ (i.e. both $arg(\zeta(1/2 + i\rho)$ and $arg(\zeta'(1/2 + i\rho))$ need to be calculable). In fact, because of Equation (6.15), only $\beta \equiv arg\,(\zeta'(1/2 + i\rho))$ needs to be calculable, along with variations of $\Gamma(1/2 + i\rho)$, $\psi(1/2 + i\rho)$ and related functions which are presumably well-known numerically. This is especially effective if Equation (6.19) is employed to represent the term $\tan(\alpha)$ in Equation (7.8).

## 8. The Volchkov equivalence

With reference to Equation (8.1) below, (where it has been opined (Moll, 2010) "Evaluating (it) might be hard"), given an explicit expression for $arg(\zeta(1/2 + i\rho))$ (see Equation 6.9), it now becomes possible to investigate the *Volchkov Criterion*, the truth of which is advertised as being equivalent to RH (Conrey & Farmer, n.d.; He, Jejjala, & Minic, 2015, Eq. 2.3; Volchkov, 1995; Sekatskii, Beltraminelli, & Merlini, 2012; Borwein, Choi, Rooney, & Weirathmueller, 2008, Section 5.2). This criterion reads: "The Riemann Hypothesis is equivalent to

$$\int_0^\infty \frac{2\,t\;\arg(\zeta(1/2 + i\,t))}{(1/4 + t^2)^2}\,dt = \pi\,(\gamma - 3)\,.''  \tag{8.1}$$

Several variations of Equation (8.1), all of which are based on a similar derivation method, are presented in these references. Consider Figure 6, a numerical comparison between the right-hand side of Equation (6.9) and $arg(\zeta(1/2 + i\rho))$ spanning the first few zeros of $\zeta(1/2 + i\rho)$ along the critical line. This comparison demonstrates that the "constant" $k$ in Equation (6.9) must be carefully chosen to achieve equality between the two representations of $arg(\zeta(1/2 + i\rho))$. In Figure 6, we note the existence of a discontinuity in $arg(\zeta(1/2 + i\rho_0))$ at $\rho = \rho_0$. In fact,

$$arg(\zeta(1/2 + i\rho_0^-)) - arg(\zeta(1/2 + i\rho_0^+)) = -\pi \tag{8.2}$$

a property that is shared (within a sign change) with any analytic function at a simple zero because the real and imaginary parts of the function both change sign as the zero is traversed in the complex $\zeta$ plane. (See Section 10). The intent now is to analytically evaluate Equation (8.1) using the right-hand side of Equation (6.9) to represent $arg(\zeta(1/2 + it))$.

Substituting Equation (6.9) into Equation (8.1) leads to three interesting integrals, the first being

$$J_1 \equiv \int_0^\infty \frac{t \; \arctan(e^{t\,\pi})}{(1/4 + t^2)^2} \, dt \tag{8.3}$$

which may be evaluated by first integrating by parts,

$$J_1 = \frac{\pi}{2} + \frac{1}{4}\pi \int_0^\infty \frac{1}{(1/4 + t^2)\cosh(t\,\pi)} \, dt \tag{8.4}$$

eventually yielding

$$J_1 = \pi/2 + \frac{\pi}{2}\,\log(2) \tag{8.5}$$

after consulting Gradshteyn and Ryzhik (1980) (Eqs. 3.522(4) and 8.370). The second required integral is

$$J_2 \equiv \int_0^\infty \frac{t^2 \; \Re\left(\int_0^1 \psi(\frac{1}{2} + i\,\rho\,t)\,d\rho\right)}{(1/4 + t^2)^2} \, dt \,. \tag{8.6}$$

Following the method of Equation (6.10), interchange the integral and sum, and after evaluating the real part of the integral, obtain

$$\Re\left(\int_0^1 \psi(\frac{1}{2} + i\,\rho\,t)\,d\rho\right) = -\gamma - \frac{\arctan(2\,t)}{t} - \left(\sum_{n=1}^\infty \frac{n\,\arctan(\frac{2\,t}{2\,n+1}) - t}{n\,t}\right) \tag{8.7}$$

Apply the outer integration to the first two of the terms in Equation (8.7), and find (courtesy of Maple)

$$\int_0^\infty \frac{\gamma\,t^2}{(\frac{1}{4} + t^2)^2} \, dt = \frac{\gamma\,\pi}{2} \tag{8.8}$$

$$\int_0^\infty \frac{t\;\arctan(2\,t)}{(\frac{1}{4} + t^2)^2} \, dt = \frac{\pi}{2} \,. \tag{8.9}$$

To evaluate the third term in Equation (8.7), interchange the convergent (grouped) sum with the outer integration and evaluate the integral (Maple), yielding

cogent ·· mathematics

$$\int_0^\infty \frac{(n \, \arctan(\frac{2t}{2n+1}) - t) \, t}{n \, (\frac{1}{4} + t^2)^2} \, dt = -\frac{\pi}{2 \, (n+1) \, n} \, . \tag{8.10}$$

Sum the resulting series (Maple), incorporate Equations (8.8 and 8.9) and eventually arrive at

$$J_2 = -\gamma \pi / 2 \, . \tag{8.11}$$

The third integral, corresponding to the term $k\pi$ in Equation (6.9), is needed to account for the area between the continuous and discontinuous functions in Figure 6. Temporarily guided by Figure 6 and Equation (8.2) which indicate that each continuous segment of $\arg(\zeta(1/2 + i\rho))$ is bounded by the $k^{th}$ zero of $\zeta(1/2 + i\rho_k)$ (but see Section 9), split the integration limits into intervals, leading (Maple) to:

$$J_3 \equiv \sum_{k=1}^\infty \int_{\rho_k}^{\rho_{k+1}} \frac{2 \, t \, k \, \pi}{(\frac{1}{4} + t^2)^2} \, dt = -16 \, \pi \sum_{k=1}^\infty \frac{k \, (\rho_k{}^2 - \rho_{k+1}{}^2)}{(4 \, \rho_k{}^2 + 1) \, (4 \, \rho_{k+1}{}^2 + 1)} \, . \tag{8.12}$$

The series in Equation (8.12) is first decomposed by partial fractions and since grouped terms (partially) cancel, the final form of the series becomes

$$J_3 = \pi \sum_{k=1}^\infty \frac{1}{\rho_k{}^2 + 1/4} \tag{8.13}$$

which can be written

$$J_3 = \frac{\pi}{2} \sum_{\rho_K} \frac{1}{|1/2 + i \, \rho_K|^2} \tag{8.14}$$

where $\rho_K$ represents all zeros of $\zeta(1/2 + i\rho_K)$ that lie on the critical line, and the overall factor $\frac{1}{2}$ has been included in recognition of the fact that the sum now includes complex conjugate values that were not included in the original sum (Equation 8.13), indicated by the use of $k \rightarrow K$. A more general form of the sum (8.14) is known (Edwards, 2001, p. 159), its value being

$$\sum_\tau \frac{1}{|\tau|^2} = 2 + \gamma - \log(4\pi) \tag{8.15}$$

where the sum over $\tau = \sigma + i\rho$ includes all zeros of $\zeta(\tau)$, including any that may not lie on the critical line. If RH is true (that is, $\Re(\tau) = \frac{1}{2}$ for all $\tau$), then the sums in Equations (8.14 and 8.15) coincide, yielding

$$J_3 = \frac{\pi}{2} \, (2 + \gamma - \log(4\pi)) \, . \tag{8.16}$$

Finally, we have the less interesting integrals to evaluate, specifically (Maple)

$$J_4 = \ln(2 \, \pi) \int_0^\infty \frac{t^2}{(\frac{1}{4} + t^2)^2} \, dt = \frac{1}{2} \, \pi \ln(2 \, \pi) \tag{8.17}$$

and

$$J_5 = -\frac{9}{4} \, \pi \int_0^\infty \frac{t}{(\frac{1}{4} + t^2)^2} \, dt = -\frac{9 \, \pi}{2} \, . \tag{8.18}$$

Putting all the parts together, that is

$$J_1 - J_2 + J_3 + J_4 + J_5 \tag{8.19}$$

gives

$$\int_0^\infty \frac{2\,t\,\arg(\zeta(\tfrac{1}{2} + i\,t))}{(1/4 + t^2)^2}\,dt = (\gamma - 3)\,\pi + T_0 \tag{8.20}$$

with

$$T_0 \equiv -\frac{1}{2}\,\pi\,\overline{\sum_\tau}\,\frac{1}{|\tau|^2} \tag{8.21}$$

where the sum over $\tau$ in $T_0$ only includes zeros of $\zeta(s)$ that do not lie on the critical line, indicated by the symbol $\overline{\sum}$. The term $T_0$ arises from the difference between the terms included in the sums appearing in Equations (8.14 and 8.15), assuming RH to be false. It is worth noting that individual terms in the sum (Equation 8.20) are positive, so there is no possibility that a sum composed of non-zero terms could itself vanish—its contribution would always be negative. With $T_0 = 0$, the claim Equation (8.1) is verified, but it does NOT prove the Riemann Hypothesis which postulates that $T_0$ vanishes; it only verifies the equivalence of Equation (8.1) and the RH as embodied in Equation (8.15) because the verification is contingent on the equality of the sums in Equations (8.14 and 8.15), which itself depends on the RH.

This unsurprising result merits further discussion and raises the question of the value of so-called "RH equivalences" to prove RH. In the first place, the original wording: "RH is equivalent to Equation (8.1)" is not very well chosen because the word "equivalent" implies a two-way correspondence (of truth); the original wording should have been: "If RH is true then Equation (8.1)" without implying the converse. This interpretation can be established by examining the derivation of Equation (8.1) which was based on a contour integration about a region where zeros of $\zeta(s)$ were presumed not to exist. Fundamentally, Equation (8.1) could have been obtained in two different ways: "If RH then Equation (8.1)" or "If not RH then Equation (8.1) plus additional terms (i.e. $T_0$)". In fact, such terms are explicitly presented, but omitted in the derivation (see Sekatskii et al., 2009, unnumbered equation terminating Section 3.1, Power functions) and a long discussion in He et al. (2015). Consequently, irrespective of which assumption is used to obtain Equation (8.1), subsequent analysis (e.g. as presented here) must be done under that same assumption, and a proof (or disproof) of RH will only emerge if either assumption, used consistently throughout, yields a contradiction (reductio ad absurdum). Otherwise, the best that can be hoped for will be a tautology. That is what has happened here.

The original derivation of Equation (8.1) was performed under the first assumption (RH is true) and, if the above analysis (Section 8) had been done with that same assumption, the $T_0$ term appearing in Equation (8.19) would not have been present. The result would have been a tautology: rhs of Equation (8.1) = rhs of Equation (8.19) (without $T_0$). Alternatively, if the original derivation of Equation (8.1) had been done under the premise "RH is not true", then additional terms, exactly equivalent to $T_0$ as it appears in Equation (8.19), would have been present in Equation (8.1) in agreement with the subsequent analysis (performed above) under that same assumption. The result would be, and is, again a tautology: rhs of Equation (8.1) (plus $T_0$) = rhs of Equation (8.19). Sans contradiction, neither of these can prove or disprove RH. All of this suggests that a proof of RH will never be obtained by inventing so-called "equivalences" that depend on hypothesizing that zeros off the critical line do, or do not, exist. However, in this case, it is a useful exercise to demonstrate the validity of Equation (6.9) and conveniently suggest some RH-independent results by examining other tacit assumptions that have been incorporated into the analysis.

First, the foregoing analysis provides a simple demonstration that the sum appearing in Equation (8.12) is convergent, for the simple reason that it equals a collection of terms that are all finite, without recourse to conditions. The transition from Equation (8.12) to Equation (8.14) involves a regrouping of terms, so it does not prove that Equation (8.14) is unconditionally convergent. A complicated proof of the convergence of Equation (8.14) can be found in many sources (e.g. Edwards, 2001). Secondly, it suggests that a careful numerical evaluation of Equation (8.1) or one of its various relatives can be used to place bounds on the smallest value of $\tau$ associated with a term of the sum appearing in Equation (8.20) (see Sekatskii et al., 2015, where an exponential weight function is appended). Thirdly, it appears to validate counting theorems that place bounds on the location of the zeros (see Section 12) by the presence of discontinuities of the argument. It will now be shown that all of the above are compromised.

### 9. But...But...But

#### 9.1. A close look at the zeros

As noted above, verification of Equation (8.1) given here, depends on Equation (8.2). In general, a discontinuity of the *arg* operator for any analytic function $h(z)$ can be associated with either of two events. The first is associated with the presence of a zero of order $n$, viz.

$$h(\rho) \approx h^{(n)}(\rho_0)(\rho - \rho_0)^n.$$  (9.1)

Specializing to the case $h(\rho) = \zeta(1/2 + i\rho)$, it is easily shown that if $\zeta(1/2 + i\rho)$ possesses a zero of order $n$ at $\rho = \rho_0$, Equation (8.2) becomes

$$arg(\zeta(1/2 + i\rho_0^-)) - arg(\zeta(1/2 + i\rho_0^+)) = -n\,\pi$$  (9.2)

and observation suggests that we take $n = 1$, consistent with the assumption that $\zeta(1/2 + i\rho)$ only possesses simple zeros.

Consider Figure 7 which shows the locus of the point $(\zeta_R, \zeta_I)$ in the complex $\zeta$ plane as it passes through a typical zero. In general, for $\rho \gtrsim 10$, because the function $\alpha(\rho)$ is monotonic with negative slope, the locus of that point must follow a clockwise path. This demonstrates a fundamental difference between $arg(\zeta(1/2 + i\rho))$ and $arg(\Gamma(1/2 + i\rho))$ whose slope is positive (see Figure 1) and whose locus thereby follows a counter-clockwise path. The question arises—Is it possible that the locus of the point $(\zeta_R, \zeta_I)$ could traverse the negative real $\zeta_R$ axis and thereby generate a discontinuity because the *arg* operation is restricted to $(-\pi, \pi)$? Such a discontinuity in the $arg(\zeta(1/2 + i\rho))$ function would be distinct from that associated with a (full) zero. This could only occur in the case of a negative imaginary half-zero, that is $\zeta_I = 0, \zeta_R < 0$ for some value of $\rho$. The LogGamma function demonstrates how such a possibility arises—see Figure 1.

Now, except at a full zero, for $\rho \gtrsim 10$, Equation (6.1) requires that the function

$$v \equiv \zeta_I' \zeta_I + \zeta_R' \zeta_R < 0.$$  (9.3)

At a discontinuity, we require $\zeta_I = 0$ corresponding to $arg(\zeta(1/2 + i\rho)) = -\pi$, and $\zeta_R' \zeta_R = (\zeta_I)' \zeta_R$ according to Equation (6.3). Also, $(\zeta_I)' < 0$ because the slope of the motion of $\zeta_I$ is expected to be negative as the locus approaches $\zeta_I = 0$ on a clockwise trajectory in the negative (lower) half-plane—see Figure 7. In that case, the product $(\zeta_I)' \zeta_R < 0$ demands that $\zeta_R > 0$, contradicting the postulate that $arg(\zeta(1/2 + i\rho)) = -\pi$ which, by definition, stipulates that $\zeta_R < 0$. This means that, since one expects (see Equation (6.3) $\zeta_R' = (\zeta_I)' < 0$ in the third quadrant, imaginary half-zeros ($\zeta_I = 0, \zeta_R \neq 0$) will only occur with $\zeta_I = 0, \zeta_R > 0$. Therefore, $arg(\zeta(1/2 + i\rho)) = 0$ and no associated discontinuity is expected to occur. Such is the perceived reality incorporated into many well-accepted results, and it is difficult to visualize how this could fail to be the case.

However, the capacity of the $\zeta$-function to surprise is boundless. Consider Figure 8, which shows information similar of that of Figure 7 except that it focuses on the region bounded by $\rho = 414...416$, containing the consecutive zeros $\rho_{212} = 415.01881..$ and $\rho_{213} = 415.45521...$. In this case, where two zeros $\rho_{212}$ and $\rho_{213}$ are in close proximity, as $\rho$ increases, the locus passes first through $\rho_{212}$, loops (clockwise) within the first quadrant, then passes through $\rho_{213}$ **without crossing the positive $\zeta_R$ axis** and thereby reverses the sign of the right-hand side of Equation (9.2). Consequently, it appears that $\zeta'_R > 0$ in the third quadrant, the previous discussion reverses, and consistently, the locus traverses the negative $\zeta_R$ axis at $\rho \approx 415.6$, and a negative imaginary half-zero arises. This generates a discontinuity in $arg(\zeta(1/2 + i\rho))$ that is not associated with a zero of $\zeta$. I shall refer to this configuration by the term "anomalous zero".

**Figure 7. Clockwise travelling locus of $\zeta(1/2 + i\rho)$ in the complex $\zeta$ plane near $\rho_6 = 37.586\ldots$ as $\rho$ increases from 37 to 40.5. The locus passes through both a full-zero and a half-zero ($\zeta_I = 0, \zeta_R > 0$), as $arg(\zeta(1/2 + i\rho))$ always decreases (i.e. increases negatively).**

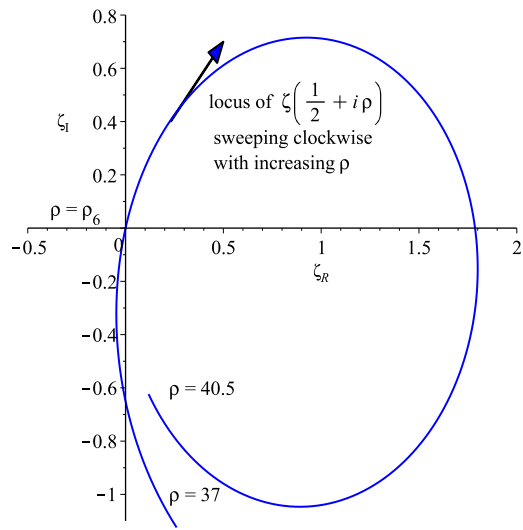

For a better understanding, consider Figure 9 where the discontinuity in $arg(\zeta(1/2 + i\rho))$ at $\rho = \rho_{213}$ shows as a **reduction** by an amount $\pi$, rather than an increase, as the locus enters the third quadrant. At $\rho \approx 415.6$, the locus crosses the negative $\zeta_R$ axis (see Figure 8), $arg(\zeta(1/2 + i\rho))$ increases by $2\pi$ and recovers to the value it would have attained if the anomalous discontinuity at $\rho_{213}$ had not occurred. For comparison, the usual situation associated with a zero of $\zeta(1/2 + i\rho)$ that results in a change in $arg(\zeta(1/2 + i\rho))$ is shown at $\rho_{212}$ corresponding to $\rho \approx 415.0$. The region that is bounded by values of $arg(\zeta(1/2 + i\rho))$ that do not adhere to their expectations is filled in solid (orange). Of note is the fact that the function $\upsilon(\rho)$ defined in Equation (9.3) is negative throughout except at the zeros, where it vanishes, as predicted, and that at the right-hand boundary of the region marked by the solid fill, the function $\zeta'_R > 0$ and $\zeta_I = 0$, consistent with the analysis given above.

Before considering the implications of these observations, it is fair to ask if there are other such anomalies. A simple search leads to a study of the pair of zeros $\rho_{126} = 279.229$ and $\rho_{127} = 282.465$ whose locus with increasing $\rho$ is shown in Figure 10. It is easier to understand the situation here by considering Figure 11 which, for clarity, is analogous to Figure 9 but shows only $arg(\zeta(1/2 + i\rho))$ and $\upsilon(\rho)$ in the region between the two zeros under scrutiny. Here, we see that the situation is the exact opposite of the previous one. In the first configuration, the two zeros $\rho_{212}$ and $\rho_{213}$ were (suspiciously) too "close" to each other relative to the point where an anomaly can occur. This separation is, in turn, regulated by the slope of $\alpha$ defined in Equation (6.6) (see Equation (13.14) below). In the second configuration, the two zeros are too "far apart" for analogous reasons, and the discontinuity rises then falls, rather than the reverse. If $\rho_{127}$ were moved 0.01 units to the left, and $\rho_{213}$ were moved 0.2 units to the right, neither anomaly would occur. However, all the location values cited here have been checked against the original published tables (Odlyzko, n.d., file zeros1) and they are consistent with those presented there, as well as those used by Mathematica with which the Maple graphics presented here have been verified. Additionally, the locations of $\rho_{212}$ and $\rho_{213}$ have been verified by an independent calculation [A. Odlyzko, private communication].

**Figure 8. Clockwise-travelling locus of $(1/2 + i\rho)$ in the complex plane near $\rho_{212} = 415.01881..$ as $\rho$ increases from 414.0 to 416.0.**

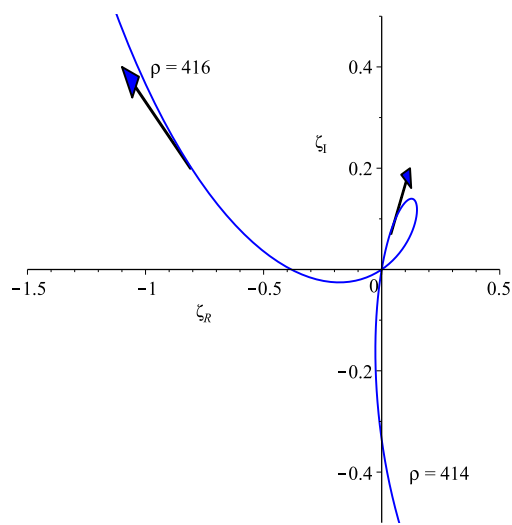

At this point, it is worthwhile to recall three formulae that locate zeros. The first is the well-known Backlund (Edwards, 2001, p. 128) counting formula

$$N(\rho) = 1 - \frac{\log(\pi)\rho}{2\pi} + \frac{\Im\left(\mathrm{Log}\Gamma\left(\frac{1}{4} + \frac{i\rho}{2}\right)\right)}{\pi} + \frac{\arg\left(\zeta\left(\frac{1}{2} + i\rho\right)\right)}{\pi} \qquad (9.4)$$

which (theoretically) increments by one at each point $\rho = \rho_k$. The dotted line $N(\rho) - 212$ in Figure 12 shows that this is not true at the anomalous point $\rho_{213}$, where it decreases by unity, before it increases by 2 at the next half-zero. Recently, França and LeClair (2013, Eq. 15) have proposed the counting function

$$N_0(\rho) = \frac{\rho}{2\pi} \ln\left(\frac{\rho}{2\pi e}\right) + \frac{7}{8} + \frac{\arg\left(\zeta\left(1/2 + i\rho\right)\right)}{\pi} \qquad (9.5)$$

based on the same (theoretical) assumption as used in Equation (9.4). The function $N_0(\rho)$ is not explicitly shown in Figure 12 because it coincides exactly with the Backlund result. In França and LeClair (2013, Remark 6), it is written "$N_0$ counts the zeros on the critical line accurately, i.e. it does not miss any zero"; Figure 12 demonstrates that there exist small ranges of the critical line where this statement is, in principle, not correct.

**Figure 9. Details of $\zeta(1/2 + i\rho)$, its argument and its derivative in the region of interest.**

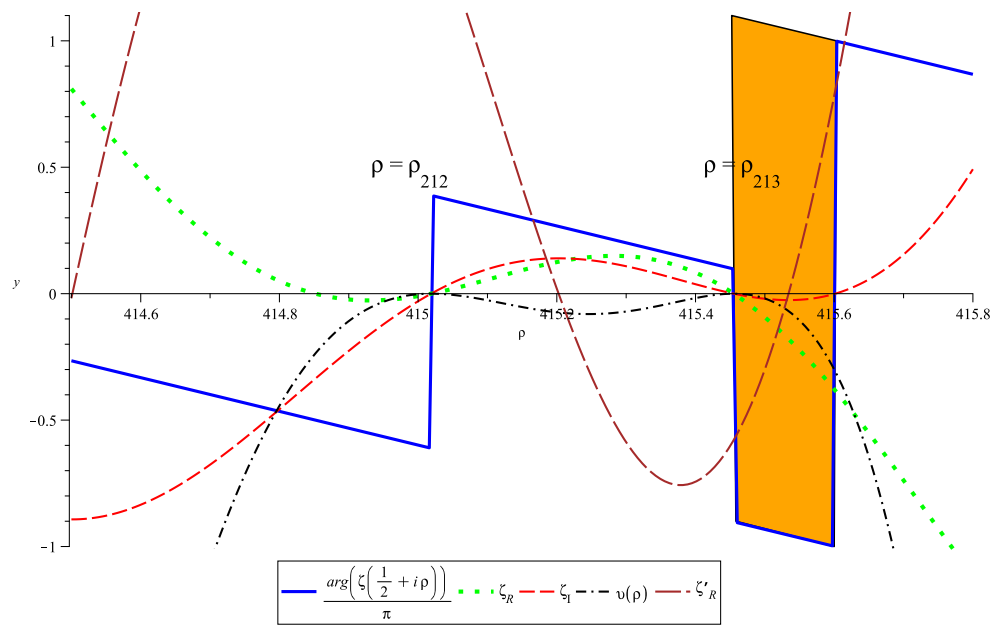

cogent •• mathematics

**Figure 10. Clockwise-travelling locus of $\zeta(1/2 + i\rho)$ in the complex $\zeta$ plane near $\rho_{127} = 282.465...$ as $\rho$ increases from 282.45 to 282.48.**

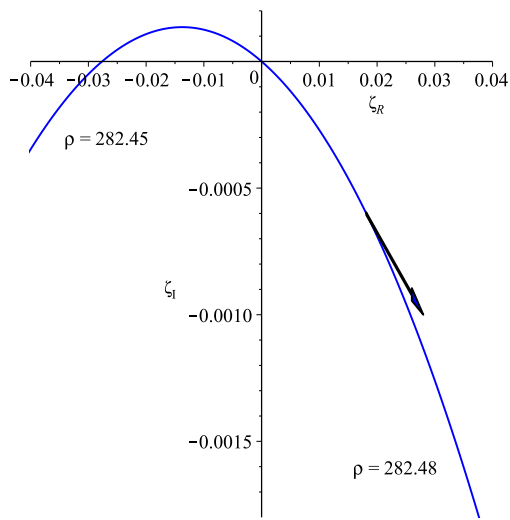

**Figure 11. Overall (left) and detail(right) of $\zeta(1/2 + i\rho)$, and $\upsilon(\rho)$ in the region of interest (see Figure 10).**

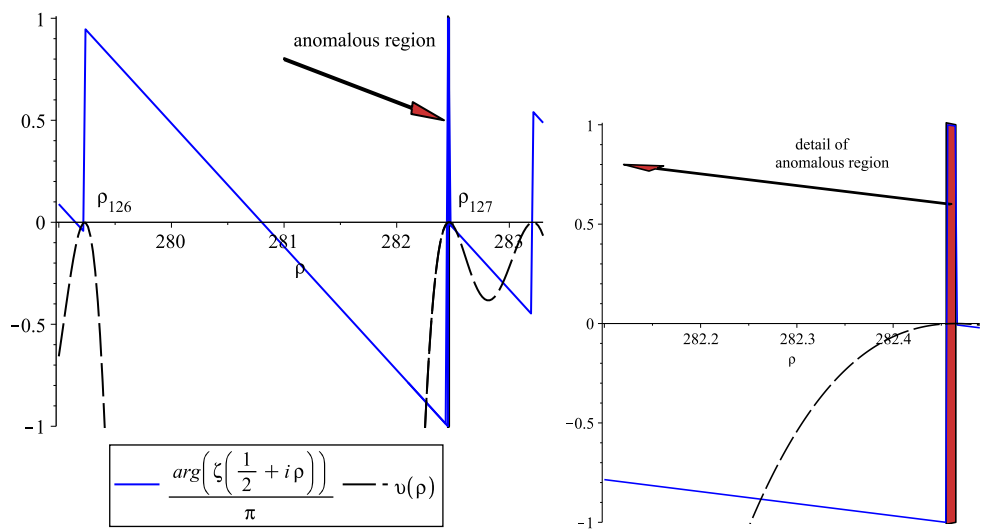

A third relevant result taken from Milgram (2011) predicts that "half-zeros corresponding to $\zeta_I = 0, \zeta_R \neq 0$ occur when $D_R = 0$" - see Appendix 1. A simple asymptotic form is easily obtained from Equation (16.1)

$$D_R \approx \frac{1}{2} - \frac{\sqrt{2}}{4}\left(\sin\left(\rho_\theta\right) + \cos\left(\rho_\theta\right)\right), \tag{9.6}$$

and Figure 12 demonstrates that the half-zero that upper bounds the anomalous zeros $\rho_{212}$ and $\rho_{213}$ is exactly where it is predicted to be according to Equation (9.6), which significantly does not depend on any numerical properties of the $\zeta$-function. The first missing term of Equation (9.6) is of order $\exp(-\pi\rho)$ and is therefore negligible when $\rho \approx 400$. The advantage of using Equation (9.6) is that numerically it only depends on simple trigonometric functions as well as the argument of $\Gamma(1/2 + i\rho)$, all of which, it is reasonably safe to assume, are numerically accurate and reliable. In this way, the calculation will not involve possible cancellation of large numbers as might occur if Equation (16.1) were to be used. Although all the anomalies observed here involve very small values of $\zeta_I$ corresponding to values of the locus travelling close to the real axis, all these suggest that the effect, being examined here, appears to be self-consistent and is unlikely to reflect numerical artefacts. In summary, we must conclude that

**Figure 12. Details of $D_R$ and $\zeta_I$ near the half-zero at $\rho = 415.601$, showing that the two coincide, and are separate from a full zero. This also shows that the counting formula Equation (9.4) fails, but only in the anomalous region.**

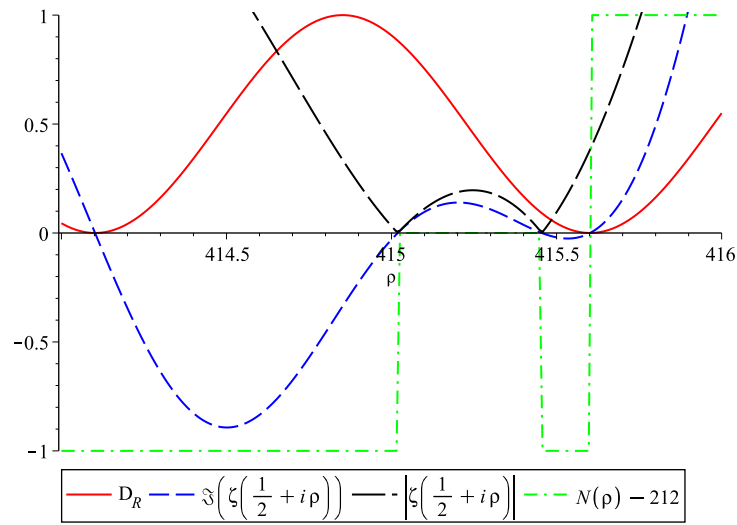

$D_R$ —— $\Im\left(\zeta\left(\frac{1}{2}+i\rho\right)\right)$ —— $\left|\zeta\left(\frac{1}{2}+i\rho\right)\right|$ —·— $N(\rho)-212$

anomalous, imaginary half-zeros of $\zeta(1/2 + i\rho)$ characterized by $\zeta_I = 0, \zeta_R < 0, \zeta_R' > 0$ do rarely, and unexpectedly, occur—a quick survey of the first 300 zeros suggests that in addition to the above, $\rho_{232-234}, \rho_{254-256}$ and $\rho_{288-290}$ may bound anomalous zeros. In addition, the Lehmer zeros (Edwards, 2001, Section 8.3) appear to share this property. As well as affecting counting formula such as Equation (9.4), the occurrence of such anomalies has implications for RH equivalences.

### 9.2. Consequences and a close look at the derivation of a ζ-function theorem

The analysis given in Section 8 depended on Equation (8.12) where the integral was split into segments, each of which was bounded by consecutive zeros $\rho_k$, justified by the assumption that each discontinuous segment of $arg(\zeta(\rho))$ was separated from the continuous function $\alpha(\rho)$ by $k\pi$, constant over that segment. Each consecutive segment was labelled by unit increments of $k$. In fact, the equivalence of Equation (8.1) and Equation (8.19), independent of the presence of the term $T_0$ as discussed, could have been taken as justification for that assumption. However, Figures 9 and 11 demonstrate that such is not the case due to anomalous zeros. With respect to Figure 9, the shaded region (orange), whose lower bound $\rho_{213}$ is not part of the area bounded by $\alpha(\rho)$ below, and $arg(\zeta(1/2 + i\rho))$ above, will nevertheless be included in the integration according to Equation (8.12). Thus, over this one segment, Equation (8.12) will over-estimate the true value of the integral by approximately

$$T_{213} \approx 2 \times (2\pi) \times (415.6 - 415.46)/415^3 \approx 2.5 \times 10^{-8} \tag{9.7}$$

and the region marked "anomalous" in Figure 11 will result in an underestimate by a numerically insignificant but non-zero result for similar reasons. Although both these anomalies are numerically insignificant, taken together with other anomalous values of $\rho_k$ that are almost certainly scattered along the critical line, one must conclude that Equation (8.1) and Equation (8.19) are analytically inconsistent, and we have a contradiction. Is this the long-sought contradiction in a RH equivalent discussed in Section 8? Unfortunately not. To understand why, it is necessary that the derivation of Equation (8.1) be revisited.

Referring to Sekatskii et al. (2012), a typical RH equivalence is obtained by integrating a function $F(z)g(z)$ around a well-chosen contour C in the complex $z$ plane, then making a convenient choice of $g(z)$. The integral is then related to the residue of enclosed singularities by the residue theorem, and it is simplified by taking one section of the contour of integration to coincide with a cut in the function $F(z)$. This is achieved by choosing $F(z) = \log(f(z))$ together with $f(z) = \zeta(z)$. The critical statement is "An appropriate choice of the branches of the logarithm function assures that the difference between the two branches of the logarithm function appearing after the integration path indents the point $(z =)X + iY$ is $2\pi i$". On the surface, this is justified by the additional requirement that "f(z) is analytic and non-zero on C". And therein lies a problem.

If $f(z)$ is analytic over a region, it does not follow that $F(z)$ will be analytic over that same region. The simplest example showing that this is true, is to consider the case $f(z) = z$. In this case, $f(z)$ is analytic, whereas $F(z) = \log(z)$ is cut. Thus, to assume that $F(z) = \log(\zeta(z))$ will retain the same analytic characteristics as $\log(z)$ and therefore integration over the same cut-structure as $\log(z)$ will yield the same discontinuity is a fairly difficult assumption to justify. And as we have seen in the previous sub-section, it is not true. There are places along the logarithm cut where the difference, defined by $F(z)$, between the two branches does not equal $2\pi i$. In fact, this occurs at precisely those points isolated previously and labelled "anomalous". Thus, if the derivation of Equation (8.1) had taken these regions into account, the statement of Equation (8.1) would have included terms precisely equal to those exemplified by Equation (9.7), thereby restoring consistency between (modified) versions of both Equation (8.1) and Equation (8.19). Thus, any RH equivalence of the Volchkov genre stated in the form Equation (8.1) without such terms is incorrect, and any numerical results based thereon must be treated with suspicion.

In a later work (Sekatskii et al., 2015), a more detailed exposition of the derivation of such equivalences is given. Using the same notation, we find written: "its (referring to $F(z)$) value jumps on $\mp2\pi l$ when we pass a point $X_1 + iY$ such that there is an $l^{th}$ order zero or pole of the function f(z) lying inside the contour (and not on the integration line) and having the ordinate Y". As has been discovered, the existence of anomalous zeros demonstrates that discontinuities can occur at points other than a zero, and this eventuality must be incorporated into the resulting proofs.

A similar caution must be issued with respect to other results where similar assumptions (propositions?) are invoked. In França and LeClair (2013, paragraph following Eq. 14), it is written "Possible discontinuities can only come from $arg(\zeta(1/2 + i\rho))/\pi$, and, in fact, it has a jump discontinuity by one whenever $\rho$ corresponds to a zero". As we have seen here, the discontinuity varies between $\pm2$ near, but not always coincident with, an anomalous zero. In Titchmarsh (1964, p. 58), we find written: "The behaviour of the function S(T) appears to be very complicated. It must have a discontinuity $k$ where T passes through the ordinate of a zero of $\zeta(s)$ of order $k$ ... Between the zeros N(T) is constant...". Again, this is not true in the vicinity of an anomalous zero.

## 10. $\zeta'(1/2 + i\rho) \neq 0$

*An immediate consequence of Equation (6.1) is:*

THEOREM 10.1   *If* $|\zeta(1/2 + i\rho)| \neq 0$, *then* $|\zeta'(1/2 + i\rho)| \neq 0$.

The proof is simple—if $\zeta_I \neq 0$ and $\zeta_R \neq 0$, then the left-side side of Equation (6.1) does not vanish, and therefore neither can the right-hand side. Except for the point $\rho = \rho_s$ (see Section 3), since the denominator (contained in the factor $f(\rho)$) does not vanish (which could cancel a numerator zero), both of $\zeta'_R$ and $\zeta'_I$ vanishing together would lead to a contradiction. Therefore, $|\zeta'(1/2 + i\rho)| \neq 0$.

This extends previously known results (Spira, 1973) onto the perforated critical line without recourse to RH.

Note 1: From Equation (2.5), $f(\rho)$ has a pole at only one point $\rho = \rho_s$, and since it is known that $\zeta(\rho_s)$ is not infinite, it must be that $\zeta'_I \zeta_I + \zeta'_R \zeta_R = 0$ at $\rho = \rho_s$ in order to cancel the pole of $f(\rho_s)$. This has been verified numerically. If $\zeta_I$ and $\zeta_R$ do not vanish, this is the only point where this can occur—otherwise, it would contradict Equation 7.8.

Note 2: It is worthwhile to recall from (Conrey, 1983a), (slightly edited): "It can be shown that RH implies that all zeros of $\xi'(s)$...have $\Re(s) = 1/2$, ... subject to simplicity (Conrey, 1983b)", and further (Conrey & Gosh, 1990), "Montgomery and Levinson proved ... $\zeta'(s)$ vanishes on $\sigma = 1/2$ only at a multiple zero of $\zeta(s)$ (hence probably never)". Titchmarsh (1964, p. 79) writes:"If Mertens' hypothesis is true, all the zeros of $\zeta(s)$ are simple".

## 11. Asymptotics

Equation (6.1) together with Equations (6.6 and 7.8) is equivalent to Equation (2.13), and taken together demonstrate that asymptotically, $|\zeta(1/2 + i\rho)|$ and $|\zeta'(1/2 + i\rho)|$ approach their limits differently as $\rho \to \infty$. That is,

$$\frac{|\zeta'(1/2 + i\rho)|}{|\zeta(1/2 + i\rho)|} \sim -\frac{\log(\rho/(2\pi))}{2\cos(\alpha - \beta)} \tag{11.1}$$

because it is known (Olver et al., 2010, Eq. 5.11.2), and easily verifiable, that

$$\Re(\psi(1/2 + i\rho)) \sim \log(\rho). \tag{11.2}$$

Within a term proportional to $\exp(-\rho\,\pi)$, Equation (11.1) is almost an equality. By way of contrast, direct differentiation of Equations (2.2 and 2.4) gives

$$\frac{|\zeta(1/2 + i\rho)|'}{|\zeta(1/2 + i\rho)|} = \frac{d}{d\rho}\ln\left(|\zeta(\tfrac{1}{2} + i\,\rho)|\right) \sim -\frac{1}{2}\ln\left(\frac{\rho}{2\pi}\right)\tan(\alpha - \beta) \tag{11.3}$$

and from Equations (11.1 and 11.3), we find a simple identity

$$\frac{|\zeta(1/2 + i\rho)|'}{|\zeta(1/2 + i\rho)'|} = \sin(\alpha - \beta). \tag{11.4}$$

Integrating Equation (11.3) between limits $T_1 \gg 0$ and $T_2 > T_1$ gives

$$\ln\left(\frac{|\zeta(\tfrac{1}{2} + i\,T_2)|}{|\zeta(\tfrac{1}{2} + i\,T_1)|}\right) \sim -\frac{1}{2}\int_{T_1}^{T_2}\ln\left(\frac{\rho}{2\pi}\right)\tan(\alpha(\rho) - \beta(\rho))\,d\rho, \tag{11.5}$$

which, lacking an explicit analytic expression for $\beta(\rho)$, "might be hard" to evaluate. Notice that if either of $T_1$ or $T_2$ coincides with a zero of $\zeta(1/2 + i\rho)$, the right-hand side of Equation (11.5) must diverge, consistent with Equation (7.8) (see Olver et al., 2010, Eq. 4.21.4), and if $T_1$ and $T_2$ enclose one or more zeros of $\zeta(1/2 + i\rho)$, the integrand acquires singularities and the integral diverges. However, Equation (11.5) has been verified by numerical integration within a few small intervals where these events do not occur.

From Equation (2.9), we establish the asymptotic limit

$$\frac{\zeta_I}{\zeta_R} \sim \frac{-\cos(\rho_\theta) - \sin(\rho_\theta) + \sqrt{2}}{\cos(\rho_\theta) - \sin(\rho_\theta)} \tag{11.6}$$

demonstrating that the real and imaginary parts of $\zeta(1/2 + i\rho)$ scale to the same order of $\rho$, subject to an (infinitely varying) modulating function given by the right-hand side of Equation (11.6). From Equation (11.1), it is also a simple matter to establish that

$$\frac{\zeta_R'}{\zeta_R} \sim -\frac{1}{2}\frac{\cos(\beta)\ln(\frac{\rho}{2\pi})}{\cos(\alpha)\cos(-\alpha + \beta)} \tag{11.7}$$

or

$$\frac{\zeta_I'}{\zeta_I} \sim -\frac{1}{2}\frac{\sin(\beta)\ln(\frac{\rho}{2\pi})}{\sin(\alpha)\cos(-\alpha + \beta)}, \tag{11.8}$$

demonstrating that corresponding components of $\zeta(1/2 + i\rho)$ and its derivative scale differently at large values of $\rho$, consistent with Equation (11.1). Numerically, both these results are excellent estimates for even reasonably small values of $\rho$ because the secondary terms that have been omitted are all of order $\exp(-\pi\rho)$.

All of these results can be understood by considering the relationship that exists between the real and imaginary components of $\zeta(1/2 + i\rho)$ and $\zeta'(1/2 + i\rho)$ presented in Equations (4.9a and 4.9b)—it is the factor $f(\rho)$ from which the scaling factor $\log(\rho/(2\pi))$ originates since all the coefficients identified in Equations (4.10 and 4.11) are $O(\rho^0)$.

Finally, with reference to Equations (6.16 and 6.17), we have, to leading orders, (and reversing the sign for convenience) the very accurate approximation,

$$\alpha_a \equiv -\alpha(\rho) \sim -\frac{\rho}{2}(1 - \log(\rho/(2\pi))) + \frac{7\pi}{8} . \tag{11.9}$$

In contrast, in Titchmarsh and Heath-Brown (1986, p. 229, Theorem 9.15), it is written: "$arg(\zeta(\sigma + i\,T)) = O(\log\,T)$ uniformly for $\sigma \geq 1/2 \ldots$". For an application, see Section 13.

## 12. Counting the zeros

Contingent on Equation (8.2) and issues raised in Section 9, we see that as $\rho$ increases past a zero, the value of $k$ in Equation (6.9) increments by one, so at any value $\rho = T$, the lower integer (floor) limit of $k$, indicated by $\lfloor \ldots \rfloor$, and computed as

$$k = \left\lfloor \frac{arg(\zeta(\frac{1}{2} + iT)) + \frac{1}{2}\int_0^T \Re(\psi(\frac{1}{2} + i\,\rho))\,d\rho - \frac{1}{2}T\ln(2\pi) + \frac{9\pi}{8} - \frac{1}{2}\,\arctan(e^{\pi T})}{\pi} \right\rfloor \tag{12.1}$$

will be equal to the number of zeros $\rho_k \leq T$. The result (Equation 12.1) is exact, but suffers from some inconsistency as an independent enumerator because it includes a (small) term that involves $\zeta$ itself (however, see França & LeClair, 2013, Remark 5). As a means of determining the location of the $k^{th}$ zero, since Equation (12.1) is exact, the floor function can be used to immediately pick out a discontinuity in $k$ that heralds the existence of a zero as $\rho$ changes, to any desired degree of numerical accuracy. As an aid to computation, Equation (12.1) can be simplified by noting that Equation (6.9) is amenable to integration by parts, yielding the alternative representation

$$\Re \int_0^T \psi(\frac{1}{2} + i\,\rho)\,d\rho = \Im(\,\log\Gamma(1/2 + i\,T)) \tag{12.2}$$

where $\log\Gamma(\ldots)$ denotes the LogGamma function (Weisstein, 2005) and, for large values of $T$, asymptotically, the right-hand side of Equation (12.2) is very well approximated (Maple). So, the first few terms of the asymptotic limit $T \to \infty$ in Equation (12.1) are

$$k = \left\lfloor \frac{1}{2}\,\frac{T\ln(T/2\,\pi)}{\pi} - \frac{T}{2\,\pi} + \frac{7}{5760\,\pi\,T^3} + \frac{1}{48\,\pi\,T} + \frac{5}{8} + 1 + \ldots \right\rfloor \tag{12.3}$$

the first two terms of which equate to the well-known result first suggested by Riemann and proven by von Mangoldt (Borwein et al., 2008, Theorem 2.9). The higher order terms can be found in Edwards (2001, Section 6.7). The last term in Equation (12.3) originates from the discontinuous term $arg(\zeta(1/2 + iT))$ in Equation (12.1), whose asymptotic limit, being unknown, perhaps unknowable, has been included in terms of its possible bounds $-\pi < arg(\zeta(1/2 + iT)) < \pi$. The traditional proof of von Mangoldt's result given at great length in Borwein et al. (2008), and discussed in Titchmarsh (1964, p. 5, Eq. 18) includes a third term $O(\log\,T)$, representing an average value of the term $arg(\zeta(1/2 + iT))$, which must vary continuously (and randomly?) between its bounds over any small range of $T$, asymptotic or not, because it is known that the number of zeros are infinite (Titchmarsh, 1964, Theorem 31). The constant term $(5/8 \pm 1)$ that is given here differs from terms given in

cogent •• mathematics

Titchmarsh ([1964](#)) and Edwards ([2001](#), Section 6.7) (i.e. $7/8 \pm \frac{1}{2}$ and 7 / 8, respectively) and is, according to Equation (6.8), an arbitrary normalization constant to begin with.

## 13. The density and distribution of zeros

In addition to the above estimates, Equation (11.9) leads to other useful approximations. Since it has been shown (see Section [9](#)) that consecutive zeros on the critical line are (almost always) bounded by half-zeros $\zeta_I = 0$, $\zeta_R > 0$, $\zeta_R' < 0$, (see Gram points (Edwards, [2001](#), p. 125)) the very accurate asymptotic estimate mod($\alpha_a$, $k\pi$) = 0 can be used to place bounds on the location of consecutive zeros since each time $\zeta_I = 0$, $\alpha_a$ will pass through a multiple of $\pi$. With $\alpha_a$ defined in Equation (11.9), Figure [13](#) illustrates this bounding procedure by plotting the normalized asymptotic counting function

$$N_a \equiv \alpha_a/\pi - \lfloor \alpha_a/\pi \rfloor \tag{13.1}$$

as a function of $\rho$ for a region of previous interest (see Section [9](#)). The absolute count value is given by $\lfloor \alpha_a/\pi \rfloor$, and the solution of the equation

$$\alpha_a = k\pi \tag{13.2}$$

which defines an upper bounding point $\rho_a(k)$ defined by $\lfloor \alpha_a/\pi \rfloor = k$, is given by

$$\rho_a(k) = (2\pi) \exp\left( W\left( \frac{-7 + 8k}{8e} \right) + 1 \right) \tag{13.3}$$

in terms of Lambert's W function (Maple). As seen in the Figure, $N_a$ bounds all zeros except the anomalous ones investigated in Section [9](#). However, as in the case of Equation (9.4), it "catches up" at the next bounding point. The accuracy of $\lfloor \alpha_a/\pi \rfloor$ as an upper bound counting function has been confirmed for the first 10,000 zeros $\rho_k$. Using the same principles, a similar result can be obtained from Equation (9.6), where $D_R = 0$ whenever

$$\theta = \frac{3\pi}{2} + 2k\pi. \tag{13.4}$$

As a comparison to Equation (13.3), França–LeClair ([2013](#)) find an expression that locates the zeros reasonable accurately, that being

$$\widetilde{\rho_k} = \frac{2\pi(k - 11/8)}{W((k - 11/8)/e)}. \tag{13.5}$$

It is interesting to compare the two approximations, by calculating the following four quantities

$$d_1 \equiv \widetilde{\rho_k} - \rho_k \tag{13.6}$$

**Figure 13.** $N_a$ as a function of $\rho$ plotted over a region of previous interest, along with $|\zeta(1/2 + i\rho)|$ to indicate the location of the zeros. The circle indicates the location of an anomalous zero.

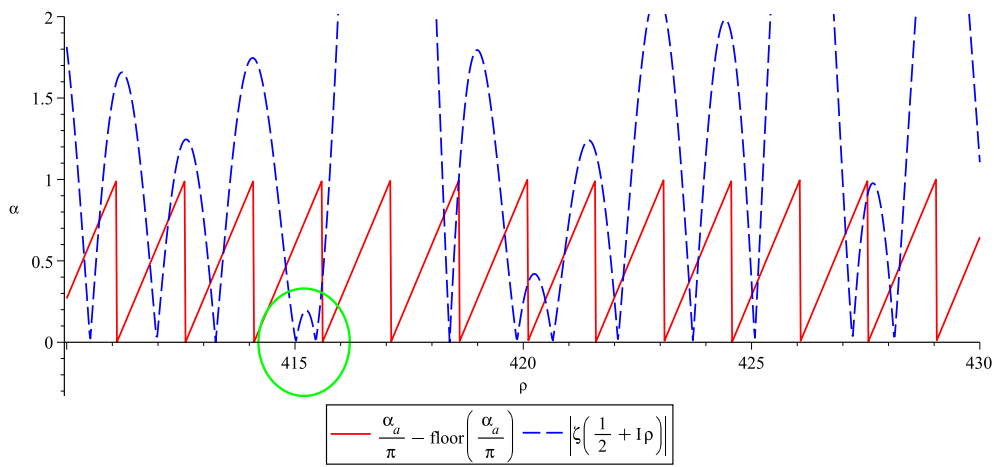

| Table 1. Statistical comparison of various measures of the distribution of zeros | | |
|---|---|---|
| | **Mean** | **Standard deviation** |
| $d_1$ | 0.00007602 | 0.2514 |
| $d_2$ | 0.4163 | 0.2514 |
| $d_3$ | -.4162 | 0.2514 |
| $d_4$ | 0.00007504 | 0.2514 |

$$d_2 \equiv \rho_a(k) - \rho_k \tag{13.7}$$

$$d_3 \equiv \rho_k - \rho_a(k-1) \tag{13.8}$$

$$d_4 \equiv (\rho_a(k) + \rho_a(k-1))/2 - \rho_k \tag{13.9}$$

each of which, respectively, measures the distance from

- the França-LeClair approximation to the known zero it is approximating;

- an upper bound to the next lower known zero;

- the lower bound up to the next higher known zero;

- the average of the upper and lower bounds to the bounded, known zero.

Such a comparison was performed using 5000 values of $\rho_k$ lying between $k = 10,000$ and $k = 15,000$. The mean and standard deviation of the different quantities, given in Table 1, show that the upper bound and lower bound are on average about 0.41 units separated from the contained zero, while both the França-LeClair approximation and an average of the upper and lower bounds are equally good. There is no observable difference in the standard deviation of any of these quantities, suggesting that each of these parameters are equally good representations of the quantity they are measuring, and that the zeros $\rho_k$ are randomly concentrated about the centre of the range bounded by $\rho_a(k)$ and $\rho_a(k-1)$.

For an estimate of the distribution of zeros at large values of $\rho$, similar to those given in França and LeClair ([2013](#)), from Equation (11.9) consider two consecutive bounding points

$$\alpha_a(\rho) = k \tag{13.10}$$

and

$$\alpha_a(\rho - \delta) = k - 1 \tag{13.11}$$

with $\delta << \rho$. Solving Equations (13.11 and 13.12) to first order in $\delta$, we find the measure of a bucket that bounds successive zeros at large values of $\rho$, that being

$$\delta \approx \frac{2\pi}{log(\rho/(2\pi))} \tag{13.12}$$

from which it is possible to estimate the density of zeros, i.e. the number of zeros $N_g$ that lie between unit intervals of $\rho$, specifically

$$N_g = 1/\delta. \tag{13.13}$$

Note that Equation (13.12) is the incarnation of Titchmarsh ([1964](#), Theorem 41) on the critical line ("the gaps between the ordinates of successive zeros of $\zeta(s)$ tend to zero"). See also Titchmarsh and Heath-Brown ([1986](#), Theorems 9.12 and 9.14). An interesting consequence of these considerations yields an estimate of the maximum possible discontinuity between successive zeros at large values of $\rho$. If, for some $\rho_a$ and $\rho$ in the range $\rho \geq \rho_a$, this maximum value, estimated by $-\delta\alpha'$, were less than $\pi$, then ALL zeros would be anomalous for $\rho \geq \rho_a$, and all existing counting theorems would fail. However, from Equations (6.6 and 13.13), we find

$$-\delta\alpha' \approx \pi \qquad (13.14)$$

so this possibility does not occur because the gap distance and the slope of the argument both vary in magnitude at the same rate (at least to first order in $\rho$).

As an example, from published tables (Odlyzko, n.d, file zeros5), at $\rho_h = 1.370919909931995308226 \times 10^{21}$ corresponding to $k = 10^{22}$, we find $N_g = 7.45$. Consistently, the tables list exactly seven zeros between $k = 10^{22} + 4$ corresponding to $\lfloor \rho_h + 1 \rfloor$ and $k = 10^{22} + 10$ corresponding to $\lfloor \rho_h + 2 \rfloor$. Similarly, there are eight listed between $k = 10^{22} + 25$ corresponding to $\lfloor \rho_h + 4 \rfloor$ and $k = 10^{22} + 32$ corresponding to $\lfloor \rho_h + 5 \rfloor$. As noted, each of these ranges corresponds to values of $\rho$ separated by one unit. The predicted gap size enveloping zeros for this range of $k$ is $\delta = 0.13416$, comparable to the average observed distance between the seven zeros which is 0.1460. In this range, the maximum possible predicted distance between zeros is $2\delta = 0.268$; the largest observed gap between the seven zeros is 0.221.

For larger values of $\rho$, the reader can estimate $N_g$ from Figure 14. For extremely large values of $\rho$, numerical evaluation of the function $W$ could be reasonably suspect since no details about the asymptotic evaluation of the $W$ function are given in Maple's FunctionAdvisor. In the case of Mathematica, we find (Weisstein, 2002) a complicated asymptotic expansion accurate to order

$$\epsilon_W = (\log(\log(\rho))/\log(\rho))^6. \qquad (15.15)$$

With $\rho \approx 10^{100}$, this suggests $\epsilon_W \approx 1.7 \times 10^{-10}$ which might affect calculations that must be done in multiple precision ($>> 100$ digit) arithmetic. This can be checked using Equations (13.2) and (11.9) which only involve the $\log$ function, presumably numerically trustworthy for large values of its argument. Define $G_{oo} = 10^{100}$. In França and LeClair (2013), with $k = G_{oo}$, Table 1 lists $\rho_{G_{oo}} = 2.8069\cdots \times 10^{98}$ to 102 digit accuracy. Substituting $\rho = \rho_{G_{oo}}$ in Equations (13.2) and (11.9), and assuming that $\rho_{G_{oo}}$ is not anomalous, we obtain $k = G_{oo} - 1$, a negligible discrepancy, given that $\rho_{G_{\infty}}$ in França and LeClair (2013) as presented is accurate to only the first digit beyond the decimal point, implying an inherent error greater than the size of the expected range boundaries (see below) at these values of $k$.

**Figure 14. Estimate of $N_g$, the number of zeros lying between $\rho$ and $\rho + 1$ as a function of $\rho$.**

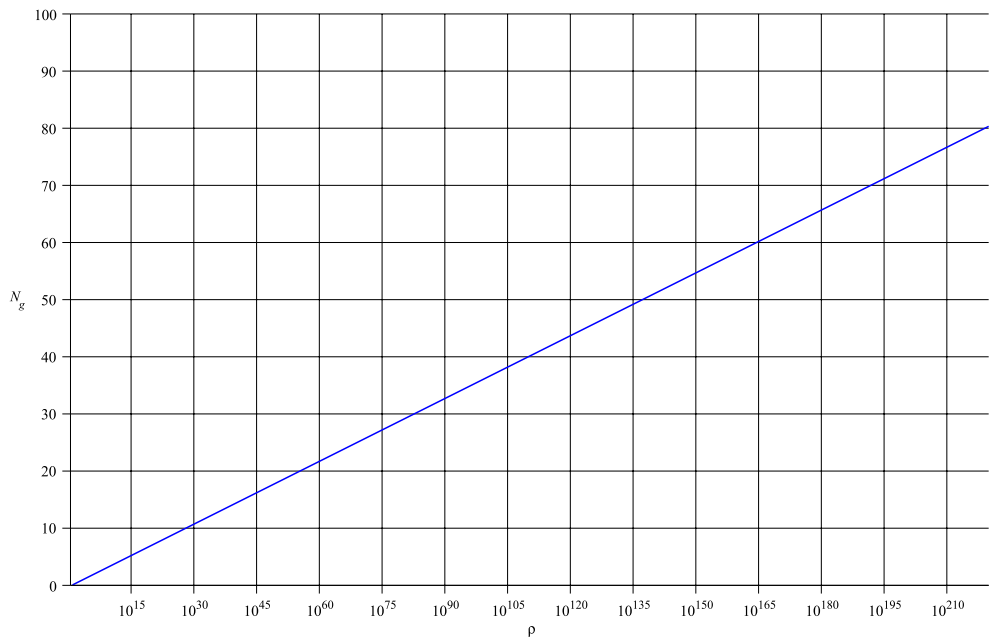

Performing the reverse test, let $k = G_{\infty}$ and $k = G_{\infty} - 1$ in Equation (13.3) to obtain the bounding points of $\rho_{G_{\infty}}$. It was found that $\rho = \rho_{G_{\infty}}$ is predicted to lie 0.013 units below the upper bound, and .014 units above the lower bound, exactly as specified in França and LeClair (2013). Furthermore, since the expected size of the distance between bounds at $\rho = \rho_{G_{\infty}}$ is 0.0279, all of this is self-consistent and appears to be acceptably accurate. Verification of the case(s) $k = 10^{200}$ (and beyond), discussed in França and LeClair (2013), is left as an exercise for the ambitious reader.

## 14. Summary

In this work, an exploratory approach has been used to obtain insights into the nature of Riemann's zeta function and its zeros in the critical strip. Specializing to the critical line, the method reproduces and expands upon a known result for locating the zeros, and yields a novel derivation of the argument of the zeta function by means of a differential equation. The discovery that a simple, singular, linear transformation exists between the real and imaginary components of $\zeta$ and $\zeta'$ on the critical line is undoubtedly significant. The analytic representation of $arg(\zeta(1/2 + i\rho))$ in turn supplies further insights, notably related to the Volchkov equivalence, a counting formula and the distribution of zeros asymptotically.

Of import, a study of the location of zeros led to the numerical discovery of anomalous zeros, and the suggestion that many previous results (e.g. Volchkov-type RH equivalents) that did not take these into account require revision. From the derived equation for $\alpha(\rho)$, it is shown, independent of RH, that $\zeta'(1/2 + i\rho)$ does not vanish on the perforated critical line. Since much has been written about the zeta function over the years, it is recognized that some of the results given here may not be new, but it is suspected that those results based on Equation (6.1) are new. At a minimum, the analyses given here do not follow usual textbook derivations and gather many disparate results in one place.

Regarding RH, it has been written (Ivić, 1985, p. 50) that "the functional equation for $\zeta(s)$ in a certain sense characterizes it completely". In fact, Titchmarsh and Heath-Brown (1986, Section 2.13) demonstrate that, with general assumptions, the functional equation defines the zeta function. Here, it has been shown that many known (and possibly new) properties of the zeta function on the critical line can be obtained by studying only the functional equation, and some results can be obtained independently of RH, in contrast to equivalent results usually cited that require RH to be true. For future investigation, a very significant advance would be accomplished by the discovery of a relationship analogous to Equation (6.6) for the argument function $\beta(\rho)$, since Equations (7.6) or or (7.8) show that the location of the non-trivial zeros on the critical line is defined by the two argument functions $\alpha$ and $\beta$ alone. The extension of Equation (13.15) to higher orders of $\rho$ would clarify how frequently anomalous zeros are expected to occur, and a proof that, in Equation (4.7), the factor enclosed in square brackets is non-zero would be useful. As well, further investigation into the consequence(s) of the existence of anomalous zeros on various accepted results is warranted. For example, the derivation of Volchkov and related equivalences and various counting theorems should be carefully and rigorously revisited.

**Supplementary material**
Supplementary material for this article can be accessed here http://dx.doi.org/10.1080/23311835.2016.1179246.

**Acknowledgements**
I am grateful to Larry Glasser, who, over the past few years, has subjected me to a barrage of novel, fascinating and wonderful insights into mathematical methods and results at a rate far faster than I could give them the attention they deserve. By a very circuitous route, the present work resulted from one such missive. Consequently, I am dedicating this paper to the memory of Larry's wife Judith Glasser. Larry, Judith and I enjoyed a memorable afternoon at a coffee shop in the border city of Cornwall, Ontario, several years ago, which led to my recent collaboration with Larry and the several papers that resulted. The author thanks Larry and Vini Anghel for a perusal of the preliminary manuscript. I also thank an anonymous referee who suggested the derivation of Equation (5.2) as it is now presented here, being orders of magnitude simpler than my original derivation.

✦ cogent •• mathematics

## Funding
The author received no direct funding for this research.

## Author details
Michael Milgram[1]
E-mail: mike@geometrics-unlimited.com
ORCID ID: http://orcid.org/0000-0002-7987-0820
[1] Geometrics Unlimited Ltd., Box 1484, Deep River, Ontario, Canada, K0J 1P0.

## Cover image
Source: Author.

## Appendix 1

## Appendix—definitions and symbols

The following symbols and definitions appear throughout this paper. Subscripts $I$ and $R$ refer to the imaginary and real components of the symbol to which they are attached. The unmodified symbols $\Gamma$, $\Upsilon$, $\psi$ and $\zeta$ mean $\Gamma(1/2 + i\rho)$, $\Gamma(1/4 + i\rho/2)$, $\psi(1/2 + i\rho)$ and $\zeta(1/2 + i\rho)$, respectively, along with their appropriate subscripts; $\psi$ is the digamma function ($\psi(z) = \frac{d}{dz}\log(\Gamma(z))$). $\zeta'_{R,I}$ refers to the real or imaginary component of the derivative of $\zeta(1/2 + i\rho)$ with respect to $\rho$, whereas $(\zeta_{R,I})'$ refers to the derivative of the real or imaginary component of $\zeta(1/2 + i\rho)$ with respect to $\rho$ (see Equations (6.3 and 6.4). All other symbols carry implicit functional dependence on $\rho$, usually omitted for brevity, except where required for clarity. The variables $n, k$ are always positive integers, and $K$ is an integer. The symbol $\rho_0$ refers to a generic value of $\rho$ corresponding to $\zeta(1/2 + i\rho_0) = 0$, whereas $\rho_k$ refers to a specific value of $\rho_0$.

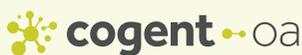

$$D_R = 1/2 - \frac{C_p \cos(\rho_\pi) + C_m \sin(\rho_\pi)}{2\sqrt{\pi}} \tag{16.1}$$

$$C_p = \Gamma_R \cosh(\frac{\pi\rho}{2}) + \Gamma_I \sinh(\frac{\pi\rho}{2}) \tag{16.2}$$

$$C_m = \Gamma_I \cosh(\frac{\pi\rho}{2}) - \Gamma_R \sinh(\frac{\pi\rho}{2}) \tag{16.3}$$

$$|\zeta'|^2 = {\zeta'_R}^2 + {\zeta'_I}^2 \tag{16.4}$$

$$\rho_\pi = \rho \log(2\pi) \tag{16.5}$$

$$\rho_\theta = \theta - \rho_\pi \tag{16.6}$$

$$\theta = arg(\Gamma(1/2 + i\rho)) \tag{16.7}$$

$$\alpha = arg(\zeta(1/2 + i\rho)) \tag{16.8}$$

$$\beta = arg(\zeta'(1/2 + i\rho)) \tag{16.9}$$

$$T_1(\rho) = -\cosh(\frac{\pi\rho}{2})\sin(2\beta + \rho_\theta)) + \sinh(\frac{\pi\rho}{2})\cos(2\beta + \rho_\theta) \tag{16.10}$$

*Cogent Mathematics* **(ISSN: 2331-1835) is published by Cogent OA, part of Taylor & Francis Group.**

**Publishing with Cogent OA ensures:**

- Immediate, universal access to your article on publication
- High visibility and discoverability via the Cogent OA website as well as Taylor & Francis Online
- Download and citation statistics for your article
- Rapid online publication
- Input from, and dialog with, expert editors and editorial boards
- Retention of full copyright of your article
- Guaranteed legacy preservation of your article
- Discounts and waivers for authors in developing regions

**Submit your manuscript to a Cogent OA journal at www.CogentOA.com**

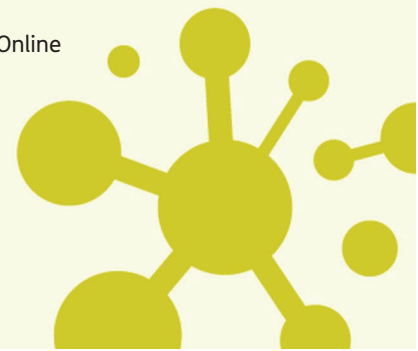